\documentclass[a4paper, 11pt, reqno]{amsart}

\usepackage{mathtools}
\usepackage[latin1]{inputenc}
\usepackage[T1]{fontenc}
\usepackage[english]{babel}
\usepackage{latexsym,amsmath,amssymb,amsfonts}
\usepackage{dsfont}
\usepackage[protrusion=false,final]{microtype}
\usepackage[dvipsnames]{xcolor}
\usepackage{graphicx}
\usepackage{float}
\usepackage{esint}
\usepackage{enumerate}
\usepackage{bbm}
\usepackage{enumitem, hyperref}
\usepackage{todonotes}
\usepackage{cleveref}
\numberwithin{equation}{section}

\usepackage{hyperref}
\hypersetup{
  colorlinks = true,
  linkcolor = blue,
  citecolor = red
}

\theoremstyle{plain}
\begingroup
\newtheorem{theorem}{Theorem}[section]
\newtheorem{lemma}[theorem]{Lemma}
\newtheorem{proposition}[theorem]{Proposition}
\newtheorem{corollary}[theorem]{Corollary}
\endgroup
\theoremstyle{definition}
\begingroup
\newtheorem{definition}[theorem]{Definition}
\newtheorem{remark}[theorem]{Remark}

\endgroup
\theoremstyle{remark}

\newcommand\be[1]{\begin{equation}\label{#1}}
\newcommand\ee{\end{equation}}

\newcommand\ep{\varepsilon}

\newcommand{\iOQ}{\int_{\Omega}\int_Q}

\newcommand{\wk}{\rightharpoonup}
\newcommand{\wkts}{\overset{2-s}{\wk}}
\newcommand{\sts}{\overset{2-s}{\to}}

\newcommand{\R}{\mathbb{R}}
\newcommand{\N}{\mathbb{N}}
\newcommand{\Z}{\mathbb{Z}}
\newcommand{\Q}{\mathbb{Q}}
\renewcommand{\S}{\mathbb{S}}
\renewcommand{\o}{\Omega}
\newcommand{\oh}{{\hat{\o}_\delta}}
\newcommand{\ohn}{{\hat{\o}_{\delta_n}}}

\newcommand{\unf}{T_{\delta_n}}
\newcommand{\bvloc}{BV_{\mathrm{loc}}}

\newcommand{\e}{\varepsilon}
\newcommand{\liminfn}{\underset{n \to \infty}{\liminf}\;}
\newcommand{\limsupn}{\underset{n \to \infty}{\limsup}\;}

\newcommand{\dhno}{\,d{\mathcal H}^{M-1}}
\newcommand{\hno}{{\mathcal H}^{M-1}}
\newcommand{\dw}{{\mathrm{d_W}}}
\newcommand{\dist}{{\mathrm{dist}}}
\newcommand{\gr}{{\mathrm{Gr}}}
\def\ca{\mathbbmss{1}}

\newcommand{\restr}{%
  \,\raisebox{-.127ex}{\reflectbox{\rotatebox[origin=br]{-90}{$\lnot$}}}\,%
}

\newcommand{\green}[1]{{\color{Green} {#1}}}
\newcommand{\blue}[1]{{\color{cyan} {#1}}}
\newcommand{\purp}[1]{{\color{violet} {#1}}}

\makeatletter
\def\namedlabel#1#2{\begingroup
    (#2)%
    \def\@currentlabel{#2}%
    \phantomsection\label{#1}\endgroup
}
\makeatother

\DeclarePairedDelimiter\floor{\lfloor}{\rfloor}

\allowdisplaybreaks

\title[Homogenization and phase separation]{Homogenization and phase separation with space dependent wells - The subcritical case}
\author[R. Cristoferi] {Riccardo Cristoferi} 
\author[I. Fonseca] {Irene Fonseca} 
\author[L. Ganedi] {Likhit Ganedi} 
\keywords{Homogenization, Phase separation, Gamma-expansion, Two-scale convergence}

\begin{document}

\maketitle

\begin{abstract}
A variational model for the interaction between homogenization and phase separation is considered. The focus is on the regime where the latter happens at a smaller scale than the former, and when the wells of the double well potential are allowed to move and to have discontinuities.
The zeroth and first order $\Gamma$-limits are identified. The topology considered for the latter is that of two-scale, since it encodes more information on the asymptotic local microstructure.
In particular, when the wells are non constant, the first order $\Gamma$-limit describes the contribution of microscopic phase separation, also in situations where there is no macroscopic phase separation.
%\red{In particular, it is shown that while there may not be macroscopic phase separation, in some cases there is a nontrivial intermediate scale where microscopic phase separation exists. In the process of proving these results, the theory of inhomogeneous Modica Mortola functionals is strengthened.}
As a corollary, the minimum of the mass constrained minimization problem is characterized, and it is shown to depend on whether or not the wells are discontinuous.
In the process of proving these results, the theory of inhomogeneous Modica Mortola functionals is strengthened.
\end{abstract}

\tableofcontents

The Modica-Mortola functional is the prototypical mathematical model for phase separation in an homogeneous material. After the initial works \cite{modica87}\cite{modica77} by Modica and Mortola that proved the conjecture by Gurtin (see \cite{Gurtin}) in the scalar case, several  variants of the functional were studied in the literature (see, for instance, \cite{sternberg88, fonseca89, KohnStern}) to prove the sharp interface limit in full generality. We refer to \cite{CriGra} for a more comprehensive overview of the gradient theory of phase separation. From the mathematical point of view, the main feature of such models is that both the potential and the wells do not depend on the spatial point, modeling homogeneity of the medium.

Modern technologies, such as temperature-responsive polymers, take advantage of engineered inclusions, or natural heterogeneities of the medium are exploited to obtain novel composite materials with specific physical properties. To model such situations by using a variational approach based on the gradient theory, the potential and the wells have to depend on the spatial point, even in a discontinuous way.

The study of phase transition in heterogeneous media with inhomogeneous conditions is a challenging mathematical problem that has recently drawn the attention of researchers.
A mathematical model for phase separation of homogeneous materials in an inhomogeneous setting was considered in the scalar case by Bouchitt\'{e} in \cite{bouchitte90}, by using techniques heavily relying on the scalar nature of the functions. A similar problem was also considered by Sternberg in \cite{sternberg88} in the two dimensional setting and with a double well required to satisfy strong regularity assumptions. The first author and Gravina in \cite{CriGra} have recently extended the above mentioned results to the vector-valued case under some strict conditions on the behavior of the double well potential around the wells.
A first result in understanding phase separation in heterogeneous media was obtained by Braides and Zeppieri in \cite{BraZep}, where the interaction between periodic microstructure and interfacial energies is studied in the scalar case in dimension one for inhomogeneous conditions. While on the one hand, the authors consider several regimes and higher order $\Gamma$ expansion, on the other hand their approach relies heavily on the explicit choices of the potential and on the wells, and on the many advantages of working in the one dimensional scalar case. In particular, the several limiting functionals identified in their work are with respect to weak-$L^2$ convergence, and the techniques used are not easily extended to the multidimensional vectorial case. In \cite{CriFonHagPop} (see also \cite{CriFonHagPop_Err}) the authors analyzed the case when the scale of the periodic microstructure and interface are of the same order, in the case of fixed wells, but without any restriction on the dimensions. Finally, we should mention that in \cite{AnsBraChi2}, the homogenization is in the singular perturbation term which leads to fundamentally different phenomenon.\\

In this paper, we consider a variational model for phase separation within a periodically heterogeneous composite material with inhomogeneous conditions, when wells may depend on the spatial variable and have discontinuities. Fixed $\e,\delta>0$, the energy can be written as
\begin{equation}\label{eq:energy_intro}
G_{\e,\delta} (u) \coloneqq \int_\o \left[\, W\left( \frac{x}{\delta}, u(x) \right)  \,+ \e^2|\nabla u(x)|^2 \,\right] dx.
\end{equation}
Here $\Omega\subset\R^N$ is an open bounded Lipschitz set, and $u\in W^{1,2}(\Omega;\R^M)$. The double well potential $W$ is $Q$-periodic in the first variable, where $Q\coloneqq(0,1)^N$, modeling a periodic structure of the material.
In the functional $G_{\e,\delta}$, the parameter $\e$ relates to the scale of the diffuse transition layer, while the scale of the periodic microstructure is $\delta$.
The main novelty of the paper is the general framework in which the asymptotic behaviour of the functional $G_{\e,\delta}$ is studied: first of all there is no restriction on the dimensions $N,M\geq 1$; moreover, for each $x\in\o$, the potential $W$ vanishes on two wells $a(x), b(x)\in\R^M$, where the $Q$-periodic functions $x\mapsto a(x)$ and $x\mapsto b(x)$ are allowed to have discontinuities.
These assumptions extend significantly those in previous works. Indeed, in \cite{CriFonHagPop} (see also \cite{CriFonHagPop_Err}) the vectorial wells were required to be fixed, and in \cite{BraZep} only the case $N=M=1$ is considered and an explicit potential $W$ and wells $a$, $b$ are used.

The core of this work is to identify the first order $\Gamma$-limit with respect to the two-scale convergence in the regime where $\e$ is negligible with respect to $\delta$, namely when the heterogeneities of the material are of a larger scale than that of the diffuse interface between different phases. The choice of working with the two-scale convergence is to maintain in the limit fine information about the asymptotic local microstructure. The second order $\Gamma$-limit, as well as other regimes, will be the content of forthcoming investigations.

Finally, we note that this regime is in particular of relevance to the biological phenomenon of lipid rafts. This is the theory that within the cell membrane there are many coexisting fluid phases consisting of various varieties of bonded lipids and disordered lipid phases. It was shown through the work of many collaborators (see \cite{Silvius2003} for a summary) that at physiological parameters, the phase separation occurs at the scale of nanometers which is inaccessible to microscopes. Furthermore, in \cite{Senugupta} it is noted that there is no macroscopic phase separation and that thermal fluctuations play a role in the formation of these nanodomains. This provides an apt setting to use the tools of homogenization to derive an effective theory for the material consisting of these nanodomains.

%%%%%%%%%%%%%%%%%%%%%%%%%%%%%%%%%%%%%%%%%%%%
%%%%%%%%%%%%%%%%%%%%%%%%%%%%%%%%%%%%%%%%%%%%
%%%%%%%%%%%%%%%%%%%%%%%%%%%%%%%%%%%%%%%%%%%%
%%%%%%%%%%%%%%%%%%%%%%%%%%%%%%%%%%%%%%%%%%%%

\subsection{Main results}

In this paper, we consider the regime $\e \ll \delta$, namely when the phase separation process happens at a lower scale than that of the heterogeneities of the material. Our main result is the integral representation of the $\Gamma$-expansion of order one of $G_{\e,\delta}$, that rigorously justifies the writing
\[
G_{\e,\delta} = G^0 + \frac{\e}{\delta} G^1 + o\left(\frac{\e}{\delta}\right).
\]
We explicitly identify the functionals $G^0$ and $G^1$. The main novelty of this manuscript is in the characterization of the scale $\frac{\e}{\delta}$ and the functional $G^1$. This functional exhibits an interaction between the periodic microstructure and phase separation. In order to maintain the information on the local microstructure, we use the notion of two-scale $\Gamma$-convergence (see \cite{cherdantsev2012two}), so that $G^1$ is defined on the space $L^1(\o;L^1(Q;\R^M))$, where $Q\coloneqq (-1/2,1/2)^N$.
% This allows for each point $x\in\o$ to retain the information on the local microstructure $L^1(Q;\R^M)$.
The first order limiting energy has the character of a bulk energy in the first variable, and of an interfacial energy in the second variable, namely it is of the form
\begin{equation}\label{eq:asymp_exp_Gamma}
G^1(u) \coloneqq \int_\Omega \widetilde{G}^1(\widetilde{u}(x,\cdot)) \, dx,
\end{equation}
where, for each $x\in\o$, $\widetilde{u}(x,\cdot)$ is the $Q$-periodic extension of the function $y\mapsto u(x,y)$, and $\widetilde{G}^1$ is the local energy of the microstructure defined, for a function $v\in \bvloc(\R^N;\R^M)$, as
\[
\widetilde{G}^1(v) \coloneqq \int_{\widetilde{Q}\cap J_v} \dw(y,v^-(y),v^+(y)) \dhno(y).
\]
Here $\widetilde{Q}\coloneqq [-1/2,1/2)^N$, $\dw$ is a degenerate geodesic distance related to the double well potential $\sqrt{W(y,\cdot)}$ (see Definition \ref{def:dw}), and $J_v$ is the jump set of the function $v$. Due to the technical nature of all of the assumptions and the definitions required to properly introduce all of the functionals above, in this section we prefer to sacrifice the rigor and to focus on commenting the peculiarities and the difficulties of the proofs. The precise assumptions are introduced and discussed in Section \ref{sec:assumptions}, while the zeroth and the first order limiting functionals are introduced in Section \ref{sec:zero} and \ref{sec:first}, respectively.

\begin{center}
\begin{figure}
\includegraphics[scale=1]{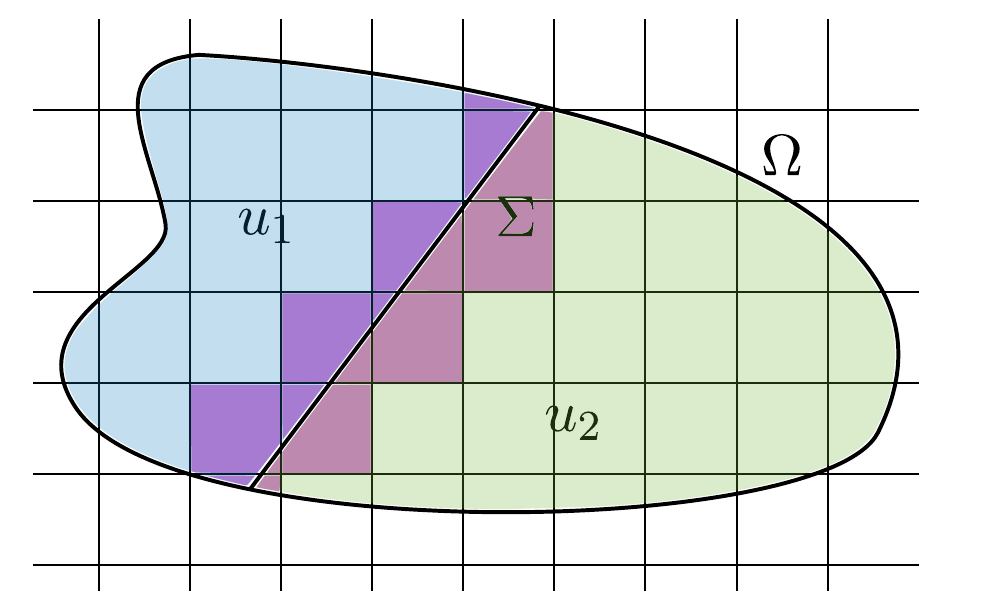}
\caption{A prototypical case showing that the energy of a recovery sequence needed to transition from a local microstructure to another is negligible in the limit. The grid is meant to represent cubes of side $\delta$, in which a different microstructure has to be approximated.}
\label{fig:Intro}
\end{figure}
\end{center}

First of all, we note that the asymptotic expansion \eqref{eq:asymp_exp_Gamma} does not depend on the rate at which $\e/\delta$ goes to zero.
Another interesting observation we gain from the form of $G^1$ is that, at first order, there is no energy penalty to pass from one local microstructure to another. To be precise, consider the situation in Figure \ref{fig:Intro}. Given a function $u:\Omega\times Q\to\R^M$ we refer to $u(x,\cdot):Q\to \R^M$ as the local microstructure at the point $x\in\o$. As expected, this function will takes values in the wells $a, b$.
Assume that $u$ is piecewise constant in the first variable, namely that it equals some $u_1:Q\to\R^M$ in the blue region, and some other function $u_2:Q\to\R^M$ in the green region. Then a recovery sequence for $u$ will have to approximate a $\delta Q$-periodic structure in the blue region, to transition between the two microstructures in the purple region, and to approximate another $\delta Q$-periodic structure in the green region. It is possible to construct the recovery sequence in such a way that the energy in each cube does not depend on the parameters $\e$ and $\delta$. Therefore, since the number of such cubes is asymptotically negligible with respect to the total number of cubes, also the energetic contribution of the recovery sequence in this region will be asymptotically negligible.
In particular, the energy of the recovery sequence will essentially be the sum of the energies needed to recover the two $Q$-periodic microstructures. This is the reason why in the functional $\widetilde{G}^1$ the jump set on $\widetilde{Q}$ is considered.

The choice of working with a potential vanishing on only two wells is based on convenience of notation: indeed, our proofs directly extends to that case of multiple wells satisfying similar assumptions as we use here.

We also note that at first order we see a local phase separation (namely in the second variable), but not a macroscopic phase separation, since this is \emph{averaged} over the entire domain. At the next order of the $\Gamma$-expansion, we expect to see a macroscopic phase separation of a similar form as the one arising from homogenization of interfaces. However, this problem will be more challenging as $\min G^1$ can be nonzero (see Corollary \ref{cor:min_energy_first}), and the structure of minimizers of the mass constrained minimization problem (which is what is most interesting for applications) might be hard to identify.

%%%%%%%%%%%%%%%%%%%%%%%%%%%%%%%%%%%%%%%%%%%%
%%%%%%%%%%%%%%%%%%%%%%%%%%%%%%%%%%%%%%%%%%%%
%%%%%%%%%%%%%%%%%%%%%%%%%%%%%%%%%%%%%%%%%%%%
%%%%%%%%%%%%%%%%%%%%%%%%%%%%%%%%%%%%%%%%%%%%

\subsection{Outline of the paper and comments on the proofs}

The paper is divided in two parts: the first is devoted to the zeroth order $\Gamma$-limit, and the second to the first order expansion.

In the realm of solid-solid phase separation, the zeroth order $\Gamma$-limit is a well-studied problem. Francfort and M\"{u}ller in \cite{FrancMull} have studied this problem in a similar framework, which was later extended by Shu in \cite{shu2000heterogeneous} to many regimes including dimension reduction. The strategy for the proof of the zeroth order $\Gamma$-limit $G^0$ we use here is, in most aspects, similar to previous work, though in the limsup inequality we employ an argument based on two-scale convergence and measurable selections.
The study of the minimization problem carried out in Corollary \ref{cor:min_F0} allows to identify the general structure of minimizers of $G^0$, which are of the form
\[
u(x) = \int_Q \mu(x,y) a(y) dy + \int_Q [1-\mu(x,y)] b(y) dy,
\]
for some $\mu\in L^2\left( \o; L^2(Q;[0,1]) \right)$. In particular, the minimum of the zeroth order asymptotic energy is zero.

Next step is to identify a class of minimizers we are interested in, and hopefully to characterize such a class by the rate of convergence to zero of the energy of a recovery sequence. Our focus will be on the class of functions that describe a geometric microstructure. Namely, those for which, for almost every $x\in\o$, the function $y\mapsto \mu(x,y)$ is a function of bounded variation taking values in $\{0,1\}$. By some heuristic computations addressed at the beginning of Section \ref{sec:first}, we get that for a function $u$ of this form, the energy of a \emph{optimal} recovery sequence is of the order $\e/\delta$. Therefore, to study the behaviour of the energy $G_{\e,\delta}$, we multiply it by $\delta/\e$, and unfold it using the two-scale unfolding operator. We study the $\Gamma$-limit of such rescaled functional which, up to a negligible error, can be written as
\begin{equation}\label{eq:G1_intro}
G^1_{\e,\delta}(u) \coloneqq \int_\o\, \left[\,\int_Q \left[\, \frac{\delta}{\e} W(y, u(x,y))
        + \frac{\e}{\delta} |\nabla_y u(x,y)|^2 \,\right] \,dy \,\right] \,dx.
\end{equation}
Note that the  nature of the limiting functional $G^1$ is clear from \eqref{eq:G1_intro}.
Compactness for sequences of uniformly bounded energy (see Lemma \ref{lem:compactness_special_class_0}) follows from an application of the Chacon biting lemma (see \cite[Lemma 2.63]{FonLeo}) together with Vitali Convergence Theorem.
The proofs of the first order $\Gamma$-limit (see Theorem \ref{thm:1_order}) use the results of \cite{CriGra} to get the \emph{inner functional} $\widetilde{G}^1$, namely that relate to the phase separation in $Q$. Since in this paper we work with more general assumptions, those cannot be directly applied, and a uniform bound on the Euclidean length of a family of geodesic problems has to be proved. Section \ref{sec:Euclidean} is entirely devoted to the proof of such bound. With this latter at our disposal, which gives the liminf inequality for the internal energy, the liminf inequality for the whole functional follows by using Fatou's lemma (see Proposition \ref{prop:liminf}).
The proof of the limsup inequality (see Proposition \ref{prop:limsup}) is based on an approximation argument. First, we consider the case where the limiting function $u\in L^1(\o;L^1(Q;\R^M))$ is piecewise constant in the first variable, namely, when
\begin{equation}\label{eq:u_intro}
u(x,y) = \sum_{i=1}^m u_i(y) \ca_{\o_i}(x),
\end{equation}
where $\o_1,\dots,\o_m$ is a polyhedral partition of $\o$, and the functions of bounded variation $u_i$'s take values on the wells. In particular, it is possible to identify each of such $u_i$'s with a set of finite perimeter $A_i\subset Q$, by setting $A_i\coloneqq\{u_i=a\}$. In this case, the recovery sequences for each of the microstructures $u_i$ provided by \cite{CriGra} are glued together in such a way that the transition between them has an asymptotically negligible energy (see Figure \ref{fig:Intro}). In order to obtain a recovery sequence for a general function $u\in L^1(\o;L^1(Q;\R^M))$, we use a density argument. This requires to being able to construct, for each $\e>0$, a function $v$ of the form \eqref{eq:u_intro} such that
\[
\| u-v \|_{L^1\times L^1} \leq \e,\quad\quad\quad
| G^1(u) - G^1(v) |\leq \varepsilon.
\]
In order to get the second inequality, as it is well known in the Calculus of Variations, the partition $\o_1,\dots,\o_m$ cannot be imposed a priori, but it has to be determined by the function $u$ itself. In particular, for measurability reasons, we need to have at our disposal a \emph{countable} family $\mathcal{C}=\{C_i\}_{i\in\N}$ of sets of finite perimeter in $Q$ such that
\[
|\widetilde{G}^1(A_i) - \widetilde{G}^1(C_k) | \leq \varepsilon
\]
for some $j\in\N$, where we naturally see the functional $\widetilde{G}^1$ as a geometric functional. The family $\mathcal{C}$ is constructed in Lemma \ref{lem:Lambda}.

The first order $\Gamma$ expansion can also be considered with respect to the weak-$L^2$ topology (see Corollary \ref{cor:main_result_x}). Moreover, the proofs we present are stable for the addition of a mass constraint to the functional (see Corollary \ref{cor:main_result_mass}). Finally, the minimization problem for the functional $G^1$ is investigated in Corollary \ref{cor:min_energy_first}.

%%%%%%%%%%%%%%%%%%%%%%%%%%%%%%%%%%%%%%%%%%%%
%%%%%%%%%%%%%%%%%%%%%%%%%%%%%%%%%%%%%%%%%%%%
%%%%%%%%%%%%%%%%%%%%%%%%%%%%%%%%%%%%%%%%%%%%
%%%%%%%%%%%%%%%%%%%%%%%%%%%%%%%%%%%%%%%%%%%%

\section{Preliminaries}

\subsection{Two-scale convergence and unfolding}
Two-scale convergence is a powerful tool, first introduced by Nguetseng \cite{Ngu} and developed further in \cite{allaire}. Later, it was separately established by Visintin \cite{Vis06,Vis07} and Cioransecu, Damlamian, and Griso \cite{CioDamGri08} to be equivalent to a topology on the product space via the use of an 'unfolding' operator. We present here some definitions and basic results obtained in the above references, which we will use in the sequel. 

We begin with the classical definitions of weak and strong two-scale convergence.

\begin{definition}\label{def:2-scale}
We say that $\{u_\delta\}_{\delta>0}\subset L^2(\Omega;\R^M)$ \emph{weakly two-scale converge to $v$} in $L^2(\Omega;L^2(Q;\R^M))$, and we write $u_\delta\wkts v$, if
\[
\lim_{n\to\infty}\int_{\Omega} u_{\delta}(x)\cdot \varphi\Big(x,\frac{x}{\delta}\Big)\,dx
=\int_{\Omega}\int_Q v(x,y)\cdot \varphi(x,y)\,dy\,dx
\]
for every $\varphi\in L^{2}(\Omega;C_{\rm per}(Q;\R^M))$. 
Here $C_{\rm per}(Q;\R^M)$ is the space of periodic continuous functions on $\R^N$ with period $Q$.
\end{definition}

Two-scale convergence encodes more information than classical weak-$L^2$ convergence. This property is highlighted in the following compactness result.

\begin{proposition}
\label{prop:2-scale-compactness}
Let $\{u_\delta\}_{\delta>0} \subset L^2(\Omega;\R^M)$ be bounded. Then, there exists $v\in L^2(\Omega;L^2(Q;\R^M))$ such that, up to the extraction of a (not relabeled) subsequence, $u_\delta\wkts v$. Additionally,
$$u_\delta\wk u\coloneqq \int_Q v(x,y)\,dy$$
weakly in $L^2(\Omega;\R^M)$.
\end{proposition}

Now we recall the  unfolding operator, with a definition that is tailored to the use that we will make of this tool.

\begin{definition}
For $\delta >0$, let
\[
\oh\coloneqq\bigcup_{z_i\in I_\delta} \left(\overline{z_i+\delta Q}\right)\cap \o,\quad\quad\quad\quad
\Lambda_\delta\coloneqq \o\setminus\oh,
\]
where $I_\delta$ is the set of points $z_i\in\delta\Z^N$ such that $\overline{z_i+\delta Q}\subset\o$.
The unfolding operator $T_{\delta}:L^2(\Omega;\R^M)\to L^2(\o;L^2(Q;\R^M))$ is defined as
\be{eq:unfolding-operator}
T_{\delta}(u)(x,y)\coloneqq\left\{
\begin{array}{ll}
u\Big(\delta\floor[\Big]{\frac{x}{\delta}}+\delta y\Big) & \text
{ for }x\in\oh,\; y\in Q, \\
& \\
a(y) & \text{ if } x\in\Lambda_\delta, y\in Q.
\end{array}
\right.
\ee
where, given an enumeration $\{z_i\}_{i\in\N}$ of $\Z^N$,
\begin{equation}\label{eq:floor}
\floor{x}\coloneqq z_i\,\quad\quad\quad i\coloneqq \min\left\{j\in\N \,:\, z_j\in \mathrm{argmin}\{ |z-x| : z\in\Z^N \}\,\right\}
\end{equation}
is the integer part of $x\in\R^N$, and $a:\Omega\to\R^M$ is the function given in \ref{W3} in Section 3.
\end{definition}

\begin{remark}
This definition of the unfolding operator is nonstandard as we make the unfolding operator nonzero in the small boundary set $\Lambda_\delta \times Q$. While this prevents our definition of the unfolding operator from being linear, it still preserves the main compactness property (see \Cref{thm:equivalent-two-scale}), and allows to simplify some algebraic arguments.
\end{remark}

\begin{theorem}
\label{thm:equivalent-two-scale}
Given $\{u_\delta\}_{\delta>0}  \subset L^2(\Omega;\R^M)$ and $v\in L^2(\Omega;L^2(Q;\R^M))$, the following conditions are equivalent:
\begin{enumerate} 
\item[(i)]$u_\delta\wkts v\quad\text{weakly two scale in }L^2(\Omega;L^2(Q;\R^M)),$
\item[(ii)] $T_{\delta}u_\delta\wk v$ weakly in $L^2(\Omega;L^2(Q;\R^M))$.
\end{enumerate}
\end{theorem}

Finally, we use the unfolding operator to define a variant of two-scale convergence that will prove useful to proving our results.

\begin{definition}
A sequence $\{u_\delta\}_{\delta>0} \subset L^1(\o;\R^M)$ is said to converge \textit{strongly two-scale in $L^1(\Omega;L^1(Q;\R^M))$} to $u \in L^1(\Omega;L^1(Q;\R^M))$ if $T_\delta u_\delta \to u$ strongly in $L^1(\Omega;L^1(Q;\R^M))$.
\end{definition}

%%%%%%%%%%%%%%%%%%%%%%%%%%%%%%%%%%%%%%%%%%%%
%%%%%%%%%%%%%%%%%%%%%%%%%%%%%%%%%%%%%%%%%%%%

\subsection{Measurable Selection}
Here, we recall the measurable selection theorem tailored to our usage that will be needed for the construction of the recovery sequences. For further reference, see \cite[Section 6.1]{FonLeo}.

\begin{definition}\label{def:multi_measurable}
A multifunction $\mathcal{F} : \Omega \to 2^Y\setminus \{\emptyset\}$, where $Y$ is a topological space, is said to be \emph{Lebesgue measurable} if for every closed set $S \subset Y$ the set
$$\mathcal{F}^{-}(S)\coloneqq  \left\{x \in \Omega\;:\; \mathcal{F}(x)\cap S \neq \emptyset\right\}$$
belongs to the Lebesgue $\sigma-$algebra on $\Omega$.
\end{definition}

The next result ensures the existence of a measurable selection (see \cite[Theorem 6.10]{FonLeo}).

\begin{theorem}\label{thm:meas_selection}
Let $Y$ be a complete separable metric space, and consider a Lebesgue measurable multifunction $\mathcal{F} : \Omega \to 2^Y\setminus \{\emptyset\}$ be with values on a closed subsets of $Y$.
Then there exists a sequence of Lebesgue measurable selections for $u_n: \Omega \to Y$, $n \in \N$, such that $\{u_n(x)\}_{n}$ is dense in $\mathcal{F}(x)$ for every $x \in \Omega$ 
\end{theorem}

%%%%%%%%%%%%%%%%%%%%%%%%%%%%%%%%%%%%%%%%%%%%
%%%%%%%%%%%%%%%%%%%%%%%%%%%%%%%%%%%%%%%%%%%%

\subsection{Sets of finite perimeter}

We recall the definition and some well known facts about sets of finite perimeter. For more details, we refer the reader to \cite{AFP, EG, M}.

\begin{definition}
Let $E\subset\R^M$ with $|E|<\infty$, and let $A\subset\R^M$ be an open set.
We say that $E$ has \emph{finite perimeter} in $A$ if
\[
P(E;A)\coloneqq\sup\left\{\, \int_E \mathrm{div}\varphi \,d x \,:\, \varphi\in C^1_c(A;\R^M)\,,\, \|\varphi\|_{L^\infty}\leq1  \,\right\}<\infty.
\]
\end{definition}

\begin{remark}\label{rem:defvar}
A set $E\subset\R^M$ is a set of finite perimeter in $A$ if and only if $\chi_E\in BV(A)$, \emph{i.e.}, the distributional derivative $D\chi_E$ is a finite vector valued Radon measure in $A$, with
\[
\int_{\R^M} \varphi \,d D\chi_E=\int_E \mathrm{div}\varphi \,d x
\]
for all $\varphi\in C^1_c(A;\R^M)$, and $|D\chi_E|(A)=P(E;A)$. In particular, the outer regularity property of Radon measures, allows to define the perimeter of a set $E\subset\R^M$ in a Borel set $D\subset\R^M$ as
\[
P(E;D)\coloneqq \inf\left\{\, P(E;A) \,:\, D\subset A, A \text{ open }  \,\right\}.
\]
\end{remark}

\begin{definition}\label{def:mun}
Let $E\subset\R^M$ be a set of finite perimeter in the open set $A\subset\R^M$. We define $\partial* E$, the \emph{reduced boundary} of $E$, as the set of points $x\in\R^M$ for which the limit
\[
\nu_E(x)\coloneqq -\lim_{r\to0}\frac{D\chi_E(x+rQ)}{|D\chi_E|(x+rQ)}
\]
exists and is such that $|\nu_E(x)|=1$.
The vector $\nu_E(x)$ is called the \emph{measure theoretic exterior normal} to $E$ at $x$.
\end{definition}

We now recall the De Giorgi's structure theorem for sets of finite perimeter.

\begin{theorem}\label{thm:DeGiorgi}
Let $E\subset\R^M$ be a set of finite perimeter in the open set $A\subset\R^M$.
Then
\begin{itemize}
\item[(i)] for all $x\in\partial* E$ the set $E_r\coloneqq \frac{E-x}{r}$ converges locally in $L^1(\R^M)$ as $r\to0$ to the halfspace orthogonal to $\nu_E(x)$ and not containing $\nu_E(x)$;
\item[(ii)] $D\chi_E=-\nu_E\,\hno\restr\partial* E$;
\item[(iii)] the reduced boundary $\partial* E$ is $\hno$-rectifiable, \emph{i.e.},
there exist Lipschitz functions $f_i:\R^{M-1}\to\R^M$, $i\in\N$, such that
\[
\partial* E=\bigcup_{i=1}^\infty f_i(K_i),
\]
where each $K_i\subset\R^{M-1}$ is a compact set.
\end{itemize}
\end{theorem}

We now recall a strong approximation result by Gromard (see \cite{DeGrom_2}, and also \cite{DeGrom_1}).

\begin{theorem}\label{thm:strong_BV}
Let $A\subset \R^M$ be an open set, and let $E\subset A$ be a set of finite perimeter in $A$.
Then, for each $\varepsilon>0$ there exist a set $F\subset A$ of finite perimeter in $A$, and a compact set $C\subset A$ such that
\begin{itemize}
\item[(i)] $\partial F\cap A$ is contained in a finite union of $C^1$ hypersurfaces;
\item[(ii)] $\| \ca_{E} - \ca_{F} \|_{BV(A)}<\varepsilon$;
\item[(iii)] $\hno(\partial F\cap A\setminus \partial* E)<\varepsilon$;
\item[(iv)] $F\subset E+B(0,\varepsilon)$, and
    $D\setminus F\subset (A\setminus E)+B(0,\varepsilon)$;
\item[(v)] $C \subset A\cap \partial* E \cap \partial F$;
\item[(vi)] $\nu_E(x) = \nu_F(x)$ for all $x\in C$;
\item[(vii)] $|D\ca_E|(D\setminus C)<\varepsilon$.
\end{itemize}
\end{theorem}

%%%%%%%%%%%%%%%%%%%%%%%%%%%%%%%%%%%%%%%%%%%%
%%%%%%%%%%%%%%%%%%%%%%%%%%%%%%%%%%%%%%%%%%%%
%%%%%%%%%%%%%%%%%%%%%%%%%%%%%%%%%%%%%%%%%%%%
%%%%%%%%%%%%%%%%%%%%%%%%%%%%%%%%%%%%%%%%%%%%

\section{Assumptions}\label{sec:assumptions}

Let $\o\subset\R^N$ be a bounded open set, and $N, M \geq 1$.
Denote by $Q\coloneqq (-1/2,1/2)^N$ the unit cube in $\R^N$, and write $\mathcal{L}^N$-a.e. point $x\in\o$ as
\begin{equation}\label{eq:xy}
x=\lfloor x \rfloor +y,
\end{equation}
where $y\in Q$, and $\lfloor x \rfloor$ is the integer part of $x\in\R^N$, defined in \eqref{eq:floor}.
Consider measurable functions $W:\o\times\R^M\to[0,+\infty)$ and $a,b:\o\to\R^M$, and pairwise disjoint open sets $E_1,\dots,E_k\subset Q$ with piecewise affine boundary and with
\[
Q=\bigcup_{i=1}^k \left( \overline{E_i}\cap Q \right),
\]
such that the following assumptions are satisfied:

\begin{enumerate}[label= (W\arabic*)]
\item \label{W1}  For all $p\in\R^M$, the function $x\mapsto W(x,p)$ is $Q$-periodic;
\item \label{W2}  For $\mathcal{L}^N$-a.e. $x\in\o$, by using the writing in \eqref{eq:xy}, it holds
\[
W(x,p) = \sum_{i=1}^k \chi_{E_i}(y) W_i(y,p),
\]
where, for each $i \in \{1, \dots, k\}$, the function $W_i:\overline{E_i}\times\R^M\to[0,\infty)$ is locally Lipschitz continuous;
\item \label{W3}  For $\mathcal{L}^N$-a.e. $x\in\o$, it holds
\[
a(x) = \sum_{i=1}^k \chi_{E_i}(y) a_i(y),\quad\quad\quad
b(x) = \sum_{i=1}^k \chi_{E_i}(y) b_i(y),
\]
where, for each $i\in\{1,\dots,k\}$, the functions $a_i, b_i:\overline{E_i}\times\R^M\to[0,\infty)$ are Lipschitz continuous. Moreover, for $\mathcal{L}^M$-a.e. $y\in\o$, it holds
\[
W(x,p)=0\quad \text{ if and only if }\quad p\in\{a(x), b(x) \};
\]
\item \label{W4} For every $i\in\{1,\dots,m\}$, and for $\mathcal{L}^M$-a.e. $y_0\in Q\setminus\{a_i=b_i\}$, there exist $\mu>0$, $R >0$, $c_1 > 0$, such that, for all $y\in B(y_0,\mu)$,
\[
\frac{1}{c_1} |p-a_i(y)|^2\leq W_i(y,p) \leq c_1 |p-a_i(y)|^2,
\]
if $|p-a_i(y)|\leq R$, and
\[
\frac{1}{c_1} |p-b_i(y)|^2\leq W_i(y,p) \leq c_1 |p-b_i(y)|^2,
\]
if $|p-b_i(y)| \leq R$;
\item \label{W5}  There exists $c_2>0$ such that, for $\mathcal{L}^N$-a.e. $x\in \o$,
\[
W(x,p)\geq \frac{1}{c_2}|p|^2,
\]
if $|p|\geq c_2 $, and
\[
W(x,p)\leq c_2(1+|p|^2),
\]
for every $p\in\R^M$.
\end{enumerate}

We would like to make several comments on the assumptions made.

\begin{remark}[Axes of periodicity]
The periodicity of the potential $W$ and of the wells $a$ and $b$ are meant to model a situation of a periodic medium, while the dependence on $x$ of all of the above functions allow to consider more general physical settings, like a  inhomogeneous one.
The choice of asking for periodicity with respect to the principal axes is not restrictive: our results hold also in the case of periodicity with respect to any basis of $\R^N$.
\end{remark}

\begin{remark}[Discontinuities of the wells]
The possible discontinuities of the potential $W$ and on the wells $a$ and $b$ make the model suitable for studying material inclusion. Our framework includes the assumptions used in the work by Braides and Zeppieri \cite{BraZep}, where $W:\R\times\R\to[0,\infty)$ is given by
\[
W(y,s)\coloneqq
\left\{
\begin{array}{ll}
\widetilde{W}(s-k) & y\in\left(0,\frac{1}{2}\right), \\
\widetilde{W}(s+k) & y\in\left(\frac{1}{2},1\right), \\
\end{array}
\right.
\]
with $\widetilde{W}(t)\coloneqq \min\{ (t-1)^2, (t+1)^2 \}$, and thus the wells are
\[
a(y) = \begin{cases}
1-k & \textup{for }y \in \left(0,\frac{1}{2}\right),\\
1+k & \textup{else},
\end{cases}, \quad b(y) =
\begin{cases}
-1-k & \textup{for }y \in \left(0,\frac{1}{2}\right).\\
-1+k & \textup{else}
\end{cases}
\]
Finally, note that there is no loss of generality in assuming that the partition $E_1,\dots,E_k$ is the same for all of the functions, and that we are also including the case where the previous functions are continuous along $\partial E_i\cap \partial E_j$.
\end{remark}

\begin{remark}[Assumptions on the wells]
We note that in the work \cite{CriGra} of the first author and Gravina, a stronger condition than \ref{W2} was assumed, i.e, that the potential $W$ is exactly \emph{quadratic} near the well. This restriction is relaxed in here, by asking only for quadratic bounds. Moreover, here we also allow wells to merge, namely we do not impose them to be well separated.
\end{remark}

\begin{remark}[On the sets $E_i$'s]
Assuming that the sets $E_i$'s to have piecewise affine boundaries is just for reader's convenience. Indeed, the only technical point where we use this assumption is in the construction of the recovery sequence for the first order $\Gamma$-limit (see Proposition \ref{prop:limsupQ}). In particular, the piecewise affine regularity of the $\partial E_i$'s allows us to apply directly the limsup inequality proved in \cite{CriGra}. For a partition with piecewise $C^1$ boundaries, a careful adaptation of the argument used to prove \cite[Proposition 4.3]{CriGra}, should give the result also in that case. Finally, if the boundaries are only Lipschitz continuous, then a Lusin type approximation with piecewise $C^1$ sets will allow to conclude.
\end{remark}

\begin{remark}[Growth at infinity of the potential]
Finally, the quadratic growth of $W$ can be generalized to any $q$-growth for $q>1$. If only the results for mass-constrained functional is of interest, then the growth can also be linear, as proved in \cite{Leo}. In all of this cases, the results will hold with the space $L^2$ substituted by the space $L^q$.
\end{remark}

\begin{remark}\label{rem:bounds}[Lower bound on the potential]
Using assumptions \ref{W2}, \ref{W3}, \ref{W4}, and \ref{W5}, it is possible to show that, for every $ r >0$, there exists $C_r>0$ such that
\[
\inf\left\{ W(x,p) \,:\, x\in \o,\, \min\{|p-a(x)|,  |p-b(x)| \}\geq r \right\} \geq C_r.
\]
\end{remark}

\begin{remark}[Extension to multiple wells]
Finally, we note that the choice of having two wells $a$ and $b$ is only for notational convenience. A similar result holds if any number of wells satisfying the above assumptions is considered.
\end{remark}

We are now in position to define the sequence of functionals that will be studied in this paper.

\begin{definition}\label{def:functional_Fn}
Let $\{\e_n\}_{n}, \{\delta_n\}_{n}$ be infinitesimal sequences such that
\[
\lim_{n\to\infty}\frac{\e_n}{\delta_n}= 0.
\]
For $n\in\N$, define the functional $G_n: L^2(\o;\R^M)\to[0,+\infty]$ as
\[
G_n (u) \coloneqq \begin{dcases} 
\int_\o \left[\, W\left( \frac{x}{\delta_n}, u(x) \right)  \,+ \e_n^2|\nabla u(x)|^2 \,\right] dx & \textup{if } u\in W^{1,2}(\o;\R^M),\\[5pt]
+\infty & \textup{else.}
\end{dcases}
\]
\end{definition}

\begin{remark}
The choice of writing the functionals by using sequences instead of using the notation  $G_{\e,\delta}$ is purely based on convenience, because when proving a $\Gamma$-convergence result, we would have had to fix some $\{\e_n\}_n, \{\delta_n\}_n$. In particular, note that, as long as they satisfy the required rate of convergence, the choice of sequences do not affect the results we present. Furthermore, the two scale results hold even when considering the subsequences indexed in $n$.
\end{remark}

%%%%%%%%%%%%%%%%%%%%%%%%%%%%%%%%%%%%%%%%%%%%
%%%%%%%%%%%%%%%%%%%%%%%%%%%%%%%%%%%%%%%%%%%%
%%%%%%%%%%%%%%%%%%%%%%%%%%%%%%%%%%%%%%%%%%%%
%%%%%%%%%%%%%%%%%%%%%%%%%%%%%%%%%%%%%%%%%%%%

\section{Zeroth-order $\Gamma$-expansion}\label{sec:zero}

This section is devoted to proving the zeroth order $\Gamma$-expansion of the functionals $G_n$ (see \Cref{def:functional_Fn}).
We start by introducing the limiting functional.

\begin{definition}
Define the functional $G^0: L^2(\o;\R^M)\to[0,+\infty]$ as
\[
G^0(u) \coloneqq \int_\o W^{\textup{hom}}(u(x)) dx,
\]
where, for $p \in \R^M$,
\[
W^{\textup{hom}}(p) \coloneqq \min\left\{\, \int_Q W^{**}(y, p+ \varphi(y))\, dy \,:\,
	\varphi\in L^2(\o;\R^M), \int_Q \varphi \,dy = 0  \,\right\}.
\]
Here, for each $y\in Q$, the function $p\mapsto W^{**}(y,p)$ is the convex envelope of the function $p\mapsto W(y,p)$.
\end{definition}

\begin{remark}\label{rmk:liminf}
By using the upper and lower bounds on $W$ (see \ref{W5}), it is easy to see that, for $p\in\R^M$, the minimization problem defining $W^{\textup{hom}}(p)$ has a minimizer $\phi$. To be precise, this occurs because $W$ does not take the value $+\infty$ and grows at least quadratically at infinity, which means we can find an affine function in the $p$-variable below $W(y,p)$. In this scenario, it is classical that the bipolar $W^{**}(y,\cdot)$, can be identified with the convex envelope of $W(y,\cdot)$ and it is known that functionals with convex integrands are weakly lower semicontinuous (see \cite[Propositions 6.31, 6.43, and Theorem 6.54]{FonLeo}).
\end{remark}

The main result of this section is the characterization of the zeroth order effective energy through $\Gamma-$convergence. A similar theorem has been proven before by Francfort and M\"uller \cite{FrancMull} in the case of solid to solid phase separation where they consider the same energy, with $u, \nabla u$ replaced by $\nabla u, \nabla^2 u$. While they use delicate approximation techniques to prove their result, we use two-scale convergence techniques to embody the spirit of this paper.

\begin{theorem}[$0^\text{th}$-order $\Gamma$-convergence]\label{thm:0_order}
Let $\{\e_n\}_{n}, \{\delta_n\}_{n}\subset (0,1)$ be infinitesimal sequences such that
\[
\lim_{n\to\infty}\frac{\e_n}{\delta_n} = 0.
\]
Let $\{u_n\}_{n}\subset W^{1,2}(\o;\R^M)$ with
\[
\sup_n G_n(u_n) < +\infty.
\]
Then, up to a subsequence (not relabeled), $u_n\wk u$ in $L^2(\o;\R^M)$ for some $u\in L^2(\o;\R^M)$ with $G^0(u)<\infty$.
Moreover, $G_n \stackrel{\Gamma}{\rightarrow} G^0$ with respect to the weak-$L^2$ convergence.
\end{theorem}

\begin{proof}
\textbf{Step 1: Compactness.} Consider a sequence $\{u_n\}_{n}\subset W^{1,2}(\o;\R^M)$ with
\[
\sup_n G_n(u_n) < +\infty.
\]
The lower bound \ref{W5} ensures pre-compactness in the weak-$L^2(\o;\R^M)$ topology. The fact that any cluster point has finite $G^0$ energy will follow from next step.\\

\textbf{Step 2: Liminf inequality.}
Let $\{u_n\}_n \in W^{1,2}(\o;\R^M)$ be such that $u_n \wk u$ in $L^2(\o;\R^M)$. Since $\{u_n\}_n$ is bounded in $L^2$, \Cref{prop:2-scale-compactness}, we may find a function $u \in L^2(\Omega;L^2(Q;\R^M))$ such that $$u_n \wkts u(x,y) \;\; \textup{and}\;\; u_n \wk u(x) \coloneqq \int_Qu(x,y) dy.$$

Assume that $\liminfn G_n(u_n) < + \infty$, as otherwise the inequality is satisfied trivially.
Firstly, we drop the gradient term, and rewrite the energy using the unfolding operator to obtain
\begin{align*}
    G_n(u_n) &\geq \int_\o W\left( \frac{x}{\delta_n}, u_n(x) \right) \, dx
    \geq \int_{\ohn} W^{**}\left( \frac{x}{\delta_n}, u_n(x) \right) \, dx, \\
    &= \int_\ohn\int_Q{\unf \left[W^{**}\left(\frac{x}{\delta_n}, u_n\right)\right]\; dydx} \\
    &= \int_\ohn\int_Q{W^{**}\left(y, \unf u_n\right)\; dydx},\\
    &= \iOQ{W^{**}\left(y, \unf u_n\right)\; dydx}.
\end{align*}
Note in the last equality we have used the fact that $W^{**}(y,a(y)) = 0$.
By \Cref{thm:equivalent-two-scale} and the definition of two-scale convergence, we know that $\unf u_n \wk u(x,y) = u(x) + v(x,y)$ for some $v \in L^2(\Omega;L^2(Q;\R^M))$ with $\int_Q v(x,y) dy = 0$.
Now the desired inequality comes directly from Remark \ref{rmk:liminf}, using the fact that the bipolar is convex and thus is weakly lower semicontinuous. We conclude that
\begin{align*}
    \liminfn G_n(u_n) &\geq \liminfn \iOQ{W^{**}\left(y, \unf u_n\right)\; dydx}, \\
    &\geq \iOQ{W^{**}\left(y, u(x)+v(x,y)\right)\; dydx} \geq G^0(u).
\end{align*}

\textbf{Step 3: Limsup inequality.} Firstly, we note that it suffices to show the limsup inequality for $u \in C^\infty_c(\Omega ; \R^M)$.

Indeed, for $u \in L^2(\Omega ; \R^M)$ we can find by density a sequence $\{u_j\}\subset C^\infty_c(\Omega ; \R^M)$ such that $u_j \to u$ strongly in $L^2$.
Furthermore, using the fact that the $\Gamma$-limsup is lower semicontinuous \cite[Proposition 6.8]{dalmaso93} with respect to weak-$L^2$, and that $G^0(u)$ is upper-semicontinuous with respect to strong $L^2$ convergence, we have
\begin{align*}
\Gamma-\limsupn G_n(u) &\leq \liminf_j \left[\Gamma-\limsupn G_n(u_j)\right] \\
&\leq \limsup_j G^0(u_j) \leq G^0(u).
\end{align*}

Fix $u \in C^\infty_c(\Omega ; \R^M)$ and $\varphi \in C^{\infty}_c(\Omega;C^{\infty}_{\text{per}}(Q;\R^M))$ with $\int_Q\varphi(x,y) dy = 0$, define $u_n^\varphi(x) \coloneqq  u(x) + \varphi(x,\frac{x}{\delta_n})$.
\begin{align*}
    \Gamma-\limsupn G_n(u) \leq \lim_{n\to\infty} G_n(u_n^\varphi) = \lim_{n\to\infty}  \int_\o \left[\, W\left( \frac{x}{\delta_n}, u_n^\varphi(x) \right)  \,+ \e_n^2|\nabla u_n^\varphi|^2 \,\right] 
\end{align*}
Note that $\sup_n |\nabla u_n^\varphi|^2 \leq \frac{C}{\delta_n^2}$. Thus, the gradient term disappears due to the assumption that $\frac{\e_n}{\delta_n}\to 0$ as $n\to\infty$. Due to the periodicity, we can apply the Riemann-Lebesgue lemma on the potential term to conclude that
$$\Gamma-\limsupn G_n(u) \leq \iOQ{W(y,u(x)+\varphi(x,y))\;dydx}$$

Let $\hat{F}^\varphi(u)$ be the functional defined in the right hand side of the inequality. Taking the lower semicontinuous envelope preserves inequalities, and we get that
\[
\Gamma-\limsupn G_n(u) \leq \mathrm{lsc}(\hat{F}^\varphi(u)).
\]

 $\frac{\e_n}{\delta_n}\to0$ as $n\to\infty$.
Let $\hat{F}^\varphi(u)$ be the functional defined in the right hand side of the inequality. Taking the lower semicontinuous envelope preserves inequalities, and we get that
\[
\Gamma-\limsupn G_n(u) \leq \mathrm{lsc}(\hat{F}^\varphi(u)).
\]
Applying a standard relaxation result (see \cite[Theorem 6.68]{FonLeo}) and taking the infimum over all $\varphi \in C^{\infty}_c(\Omega;C^{\infty}_{\text{per}}(Q;\R^M))$ with $\int_Q\varphi(x,y) dy = 0$, we obtain the inequality:
$$
\Gamma-\limsupn G_n(u) \leq \inf_\varphi{\left\{\iOQ{W^{**}(y,u(x)+\varphi(x,y))\;dydx}\right\}} =: \hat{G}^0(u).
$$

In order to finish the proof, we will show that $\hat{G}^0(u) = G^0(u)$.
Firstly, we note that by density it holds
$$
\hat{G}^0(u) = \inf\left\{\,
	\varphi\in L^2(\o;L^2(Q;\R^M)), \int_Q \varphi(x,y) \,dy = 0  \,\right\}.
	$$
It is easy to see that $\hat{G}^0(u) \geq  G^0(u)$. Next, we will show the other inequality.\\
Note that by Remark \ref{rmk:liminf}, for every $x \in \Omega$, we can find minimizers $\varphi^x$ such that:
$$W^{\textup{hom}}(u(x)) = \int_Q W^{**}(y, u(x)+ \varphi^x(y))\, dy.$$
We use \Cref{thm:meas_selection} in order to extract a measurable selection. By coercivity of $W^{**}$ and since $\|u\|_\infty<+\infty$, there exists $R_0>0$ such that the image of any minimizer $\varphi^x$ must be contained in the closed ball $B \coloneqq  \overline{B(0,R_0)} \subset L^2(Q;\R^M)$. We consider $Y \coloneqq  B$ equipped with the weak topology, which is metrizable since $L^2(Q;\R^M)$ is reflexive. In particular, it is separable and complete as well, and so satisfies the conditions in \Cref{thm:meas_selection}. Define the multifunction
$$\mathcal{F}(x)\coloneqq  \left\{\varphi \in Y :\int_Q\varphi \;dy = 0\;\textup{and}\; \varphi \textup{ attains the minimum } W^{\textup{hom}}(u(x)) \right \}.$$
As noted before, this is nonempty for every $x \in \Omega$.

Furthermore, we claim that for every $x \in \Omega$, $\mathcal{F}(x)$ is a closed subset of $Y$ under the weak topology. Indeed, we just need to show it is sequentially closed, so take $\{\varphi_n\}_n\subset \mathcal{F}(x)$ and suppose $ \varphi_n \wk \varphi$. Note that the zero average condition passes to the limit, and by sequential weak lower semicontinuity of the integrand, we have
\begin{align*}
    W^{\textup{hom}}(u(x)) &\leq \int_Q W^{**}(y, u(x)+ \varphi(y))\, dy\\ 
    &\leq \liminfn\int_Q W^{**}(y, u(x)+ \varphi_n(y)) = W^{\textup{hom}}(u(x)).
\end{align*}
and we conclude that $\varphi \in \mathcal{F}(x)$. The last condition needed to be checked in order to apply \Cref{thm:meas_selection} is that the multifunction is Lebesgue measurable. Take $S$ which is closed in $Y$ under the weak topology. We will show that $\mathcal{F}^-(S)$, as given in \Cref{def:multi_measurable}, is also closed. Again it suffices to check sequential closure, so we take $\{x_n\}_n \subset \mathcal{F}^-(S)$ such that $x_n \to x$ in $\o$. By definition of $\mathcal{F}^-(S)$, we can find for each $x_n$ a corresponding $\varphi_n \in \mathcal{F}(x_n) \cap S$. The sequence $\{\varphi_n\}_n$ is uniformly bounded in $L^2(Q;\R^M)$ by definition of $Y$, and so we can find a weakly converging subsequence to some limit $\varphi$. Furthermore, as $S$ is closed with respect to weak convergence, we must have $\varphi \in S$ and the zero average condition is preserved. We pass to that subsequence in both $\{x_n\}_n$ and $\{\varphi_n\}_n$, without relabeling. As $x_n \to x$ in the sense of $\R^N$, by continuity of $u$ we have $u(x_n)\to u(x)$. Once again, we can apply the sequential weak lower semicontinuity of the integrand and the upper semicontinuity of $W^{\textup{hom}}$, to get
\begin{align*}
    W^{\textup{hom}}(u(x)) &\leq \int_Q W^{**}(y, u(x)+ \varphi(y))\, dy \\
    &\leq \liminfn\int_Q W^{**}(y, u(x_n)+ \varphi_n(y)) \\
    &\leq  \limsupn W^{\textup{hom}}(u(x_n))
    \leq W^{\textup{hom}}(u(x)).
\end{align*}
Thus, all inequalities are actually equalities and we have by definition, $\varphi \in \mathcal{F}(x) \cap S$. This means that $x \in \mathcal{F}^(S)$, and the set is closed. Since the set is closed, it is a Borel set, which is contained in the Lebesgue $\sigma-$algebra. This proves that the multifunction is Lebesgue measurable, and so we have checked all the hypotheses of Theorem \ref{thm:meas_selection}. It allows us to find a measurable function $v: \o \to Y$ such that $v(x) \in \mathcal{F}(x)$.  
Furthermore, since $v(x) \in Y = \overline{B(0,R_0)} \subset L^2(Q;\R^M)$, we have that $\|v(x)\|_{L^2(Q;\R^M)} \leq R_0$. In particular, as $\Omega$ is a bounded set in $\R^N$, we have that $$\int_\o{\|v(x)\|^2_{L^2(Q;\R^M)}\;dx} < + \infty.$$

Therefore, we have that $v \in L^2(\o;L^2(Q;\R^M))$ and by definition of $\mathcal{F}(x)$, we have $\int_Q{v(x,y)dy} = 0$. Thus, $v$ is an admissible competitor for the infimum in $\hat{G}^0(u)$. We deduce that
$$\hat{G}^0(u) \leq \iOQ{W^{**}(y,u(x)+v(x,y))\;dydx} = G^0(u),$$
and this concludes the proof.
\end{proof}

It is also possible to get the explicit value of the minimum of the limiting functional $G^0$, as well as a characterization of the set of its minimizers.

\begin{corollary}[Minimizers of $G^0$]\label{cor:min_F0}
It holds that
\[
\min\{ G^0(u) \,:\, u\in L^2(\o;\R^M) \} = 0.
\]
Furthermore, $u\in L^2(\o;\R^M)$ is such that $G^0(u)=0$ if and only if
\begin{equation}
\label{eq:shape_of_p}
u(x) = \int_Q \mu(x,y) a(y) dy + \int_Q [1-\mu(x,y)] b(y) dy,
\end{equation}
where $\mu\in L^2(\o;L^\infty(Q;[0,1]))$.
\end{corollary}

\begin{proof}
\textbf{Step 1: Minima of Convex Envelope.} 
Note that $W^{**}(y, p+ \phi(y))=0$ if and only if
\begin{equation}\label{eq:def_varphi}
p+\varphi(y) = f(y) a(y) + [1-f(y)] b(y)
\end{equation}
for some $f(y)\in[0,1]$. This is due to the fact that $W$ is only zero at $a,b$ which is the minimum (see \ref{W3}). \\

\textbf{Step 2: Sufficiency.} 
First, suppose
\[
u(x)\coloneqq   \int_Q \mu(x,y) a(y) dy + \int_Q [1-\mu(x,y)] b(y) dy
\]
for some $\mu\in L^2(\o;L^\infty(Q;[0,1]))$. Consider
\[
\phi^x(y) \coloneqq  \mu(x,y) a(y) + [1-\mu(x,y)] b(y) - u(x).
\]
Note that $\phi^x \in L^2(Q;\R^M)$ with $\int_Q \phi^x(y) \,dy = 0$. 
Thus, it is admissible competitor in the minimization problem defining $W^{\textup{hom}}$. Furthermore, using convexity, we can deduce that
\begin{align*}
W^{\textup{hom}}(u(x)) &\leq \int_Q W^{**}\Big(\mu(x,y) a(y)+[1-\mu(x,y)] b(y)\Big) \;dy \\
&\leq\int_Q \mu(x,y) W^{**}(y,a(y))+[1-\mu(x,y)] W^{**}(y,b(y))\;dy = 0,
\end{align*}
where in the last inequality we used the fact that, for every $y\in Q$,
\[
W^{**}(y,a(y))=W^{**}(y,b(y))=0.
\]
Since $W\geq0$ and, in turn, $W^{\textup{hom}}\geq0$, we conclude that $W^{\textup{hom}}(u(x)) = 0$ for $\mathcal{L}^N$-a.e. $x \in \o$.
Thus, we obtain $G^0(u) = 0$.\\

\textbf{Step 3: Necessity.} 
Let $u\in L^2(\o;\R^M)$ be such that $G^0(u) = 0$. In the proof of the limsup inequality above, we showed that by a measurable selection, we can find $\varphi \in L^2(\o;L^2(Q;\R^M))$ such that $\int_Q{\varphi(x,y)dy} = 0$ and 
$$0 = G^0(u)=\hat{G}^0(u) = \int_\o\int_Q W^{**}(y, u(x)+ \varphi(x,y))\, dydx.$$
In particular, for $\mathcal{L}^N$-a.e. $y\in Q$ and $\mathcal{L}^N$-a.e. $x\in \o$, we must have that $W^{**}(y, u(x)+ \varphi(x,y)) = 0$.
Since $a,b$, and $\varphi$ are measurable, we can find $\mu \in L^2(\o;L^\infty(Q;[0,1]))$ such that
\[
u(x) + \varphi(x,y) = \mu(x,y) a(y) + [1-\mu(x,y)] b(y).
\]
Integrating \eqref{eq:def_varphi} in $Q$ and using the fact that $\varphi$ has zero average, we get \eqref{eq:shape_of_p}.
\end{proof}

Finally, as it is well-known in the general context of this work, adding a mass constraint to the problem does not require to change significantly the proof of Theorem \ref{thm:0_order}.

\begin{definition}
Let $\overline{m}\in\R^M$.
For $n\in\N$, define $\mathcal{G}_n: L^2(\o;\R^M)\to[0,+\infty]$ as
\[
\mathcal{G}_n(u;\overline{m})\coloneqq
\begin{dcases}
G_n(u) & \text{ if } u\in W^{1,2}(\o;\R^M) \text{ with } \int_\o u\, dx = \overline{m},\\
+\infty & \text{ else},
\end{dcases}
\]
and $\mathcal{G}^0: L^2(\o;\R^M)\to[0,+\infty]$ as
\[
\mathcal{G}^0(u;\overline{m})\coloneqq
\begin{dcases}
G^0(u) & \text{ if } \int_\o u \, dx = \overline{m},\\
+\infty & \text{ else}.
\end{dcases}
\]
\end{definition}

The analogous of Theorem \ref{thm:0_order} and of Corollary \ref{cor:min_F0} hold also for the mass constrained functional. The small changes needed in the proof are classical, and therefore we will not report them here (see, e.g., \cite{fonseca89}\cite{CriGra})

\begin{theorem}
Fix $\overline{m}\in\R^M$. Let $\{u_n\}_n\subset W^{1,2}(\o;\R^M)$ be such that
\[
\sup_n \mathcal{G}_n(u_n;\overline{m}) < +\infty.
\]
Then, up to a subsequence (not relabeled), $u_n\wk u$ in weak-$L^2(\o;\R^M)$ for some $u\in L^2(\o;\R^M)$ with $\mathcal{G}^0(u;\overline{m})<\infty$.
Moreover, $\mathcal{G}_n(\cdot;\overline{m}) \stackrel{\Gamma}{\rightarrow} \mathcal{G}^0(\cdot;\overline{m})$
with respect to the weak-$L^2$ convergence.
Finally,
\[
\min\{ \mathcal{G}^0(u;\overline{m}) \,:\, u\in L^2(\o;\R^M) \} = 0,
\]
and $u\in L^2(\o;\R^M)$ is such that $\mathcal{G}^0(u;\overline{m})=0$ if and only if
\[
u(x) = \int_Q \mu(x,y) a(y) dy + \int_Q [1-\mu(x,y)] b(y) dy,
\]
where $\mu\in L^2(\o\times Q;[0,1)])$ with $\int_\Omega u\,dx = \overline{m}$.
\end{theorem}

%%%%%%%%%%%%%%%%%%%%%%%%%%%%%%%%%%%%%%%%%%%%
%%%%%%%%%%%%%%%%%%%%%%%%%%%%%%%%%%%%%%%%%%%%
%%%%%%%%%%%%%%%%%%%%%%%%%%%%%%%%%%%%%%%%%%%%
%%%%%%%%%%%%%%%%%%%%%%%%%%%%%%%%%%%%%%%%%%%%

\section{First-order $\Gamma$-expansion}\label{sec:first}

In view of Corollary \ref{cor:min_F0}, we know that minimizers of $G^0$ are of the form
\begin{equation}\label{eq:general_u}
u(x) = \int_Q \mu(x,y) a(y) dy + \int_Q [1-\mu(x,y)] b(y) dy,
\end{equation}
for some $\mu\in L^2\left( \o; L^2(Q;[0,1]) \right)$.
We would like to study the behaviour of the sequence of functionals $G_n$ close to the subclass $\mathcal{R}$ of functions $u$ as in \eqref{eq:general_u} such that
\begin{equation}\label{eq:good_minimizers}
\mu(x,y) = \ca_{A(x)}(y)a(y) + \left( 1-\ca_{A(x)}(y) \right) b(y),
\end{equation}
where, for a.e. $x\in\Omega$, the set $A(x)\subset Q$ has finite perimeter.
The class $\mathcal{R}$ corresponds to \emph{geometric} microstructures, and we will see that this is the only class for which \emph{equipartition of surface energy} holds.
We proceed as follows: in next section we identify the proper scaling $\omega(\e)$ of the energy by using heuristic arguments, while rigorous arguments will be employed in the following sections to prove the $\Gamma$-expansion result.

%%%%%%%%%%%%%%%%%%%%%%%%%%%%%%%%%%%%%%%%%%%%
%%%%%%%%%%%%%%%%%%%%%%%%%%%%%%%%%%%%%%%%%%%%
%%%%%%%%%%%%%%%%%%%%%%%%%%%%%%%%%%%%%%%%%%%%
%%%%%%%%%%%%%%%%%%%%%%%%%%%%%%%%%%%%%%%%%%%%

\subsection{Heuristics for the scaling analysis}
Let $u_0\in L^1(\o;L^1(Q;\R^M))$, and let $\{u_n\}_n\subset W^{1,2}(\o;\R^M)$ be such that $u_n\wkts u_0$. By using a change of variable, and neglecting the contribution of cells that intersect $\partial\Omega$, we have that
\begin{align}\label{eq:heuristics}
G_n(u_n) &= \sum_{z\in \delta_n \Z^N} \int_{(z+\delta_n Q)\cap\Omega} \left[\,
    W\left(\frac{x}{\delta_n}, u_n(x)\right) + \e_n^2 |\nabla u_n(x)|^2 \,\right]dx \nonumber \\
    &\sim \sum_{z\in \delta_n \Z^N} \int_Q \left[\, W(y,T_{\delta_n} u_n(x,y))
	    + \left(\frac{\e_n}{\delta_n}\right)^2 |\nabla_y T_{\delta_n} u_n(x,y)|^2  \,\right] dy.
\end{align}
We will focus on the behaviour of the energy in each cube.
Fix $x\in\Omega$, $n\in\N$, and consider the function $\widetilde{u}_n:Q\to\R^M$ defined as $\widetilde{u}_n(x,y)\coloneqq \unf u_n(x,y)$.
It holds $\widetilde{u}_n\wk u_0(x,\cdot)$.

We first assume that $u_0(x,y)\in\{a(y), b(y)\}$ almost everywhere, and identify the scaling $\omega_n$ of \eqref{eq:heuristics} in each cube. Then we show that, if for a generic $u_0\in L^1(\o;L^1(Q;\R^M))$ and $\{u_n\}_n\subset W^{1,2}(\o;\R^M)$ with $u_n\wkts u_0$ the sequence $\{G_n(u_n)\}_n$ behaves like the scaling $\omega_n$, then the limiting function $u_0$ must be such that $u_0(x,y)\in\{a(y), b(y)\}$ almost everywhere.

\begin{figure}
\includegraphics[scale=0.4]{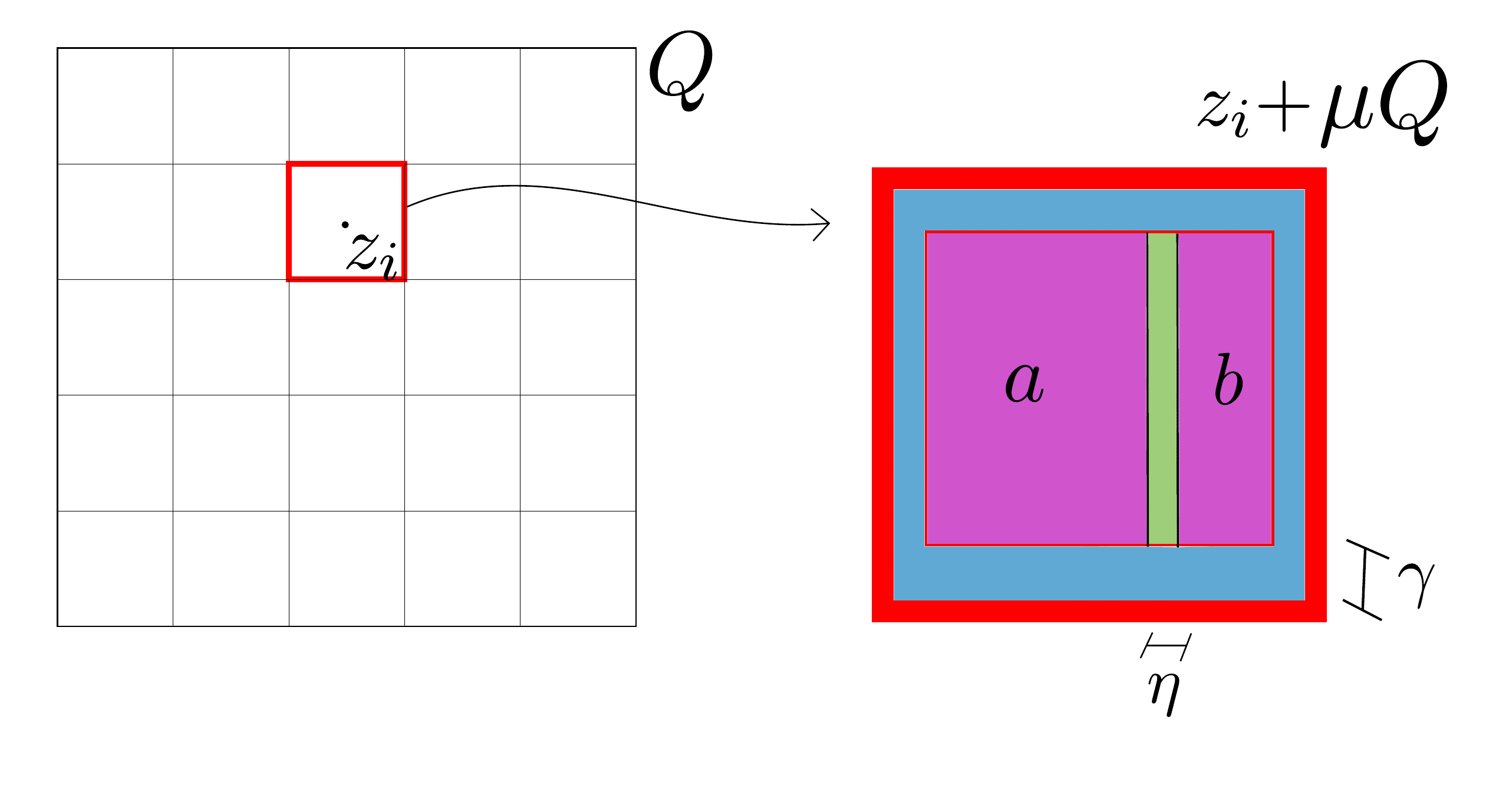}
\caption{The construction in each cube $z_i+\mu Q$: in the pink region we set $\widetilde{u}_n(y)$ to be $a(y)$, and in the purple one to be $b(y)$.
We then have the interface (colored in green) of thickness $\eta$, and the cut-off region (colored in blue) of thickness $\eta$.}
\label{fig:heuristics}
\end{figure}

Let $\mu_n>0$, to be chosen later, and subdivide the cube $Q$ into smaller cubes $\mu_n Q$. In each of these little cubes we perform the construction detailed in Figure \ref{fig:heuristics}. Namely, for $\gamma_n, \eta_n>0$ with
\begin{equation}\label{eq:constraint_eta_gamma}
\eta_n \ll \mu_n,\quad\quad\quad\quad\quad \gamma_n\ll \mu_n,
\end{equation}
the function $\tilde{u}_n(y)$ is either $a(y)$ or $b(y)$ in most of the cube, with $\eta_n$ being the thickness of the interface between $a$ and $b$, and $\gamma_n$ being the thickness of the cut-off region. Both are needed to ensure that $u_n(y) \in W^{1,2}(Q;\R^M)$.

We now evaluate the order of the energy of such a configuration. We have
\begin{align}\label{eq:scaling}
\frac{1}{\mu_n^N}& \int_{\mu_n Q} \left[\, W(y, \widetilde{u}_n(y))
	+\left(\frac{\e_n}{\delta_n}\right)^2|\nabla_n \widetilde{u}_n(y)|^2  \,\right] dy \nonumber \\
&=\frac{1}{\mu_n^N}\Bigg\{
	\purp{\left[\, (0+\left(\frac{\e_n}{\delta_n}\right)^2)(\mu_n^N-\eta_n\mu_n^{N-1}
	- \gamma_n\mu_n^{N-1}) \,\right]} \nonumber \\
&\hspace{0.5cm}+ \green{\left[\,(1 + \left(\frac{\e_n}{\delta_n}\right)^2\frac{1}{\eta_n^2})(\eta_n\mu_n^{N-1}) \,\right]}  + \blue{\left[\,(1 + \left(\frac{\e_n}{\delta_n}\right)^2\frac{1}{\gamma_n^2})(\gamma_n\mu_n^{N-1}) \,\right]} \Biggr\} \nonumber \\
&\sim \purp{\left[\,\left(\frac{\e_n}{\delta_n}\right)^2\,\right]}
+ \green{\left[\,\frac{\eta_n}{\mu_n} + \left(\frac{\e_n}{\delta_n}\right)^2\frac{1}{\eta_n\mu_n}\,\right]}
+ \blue{\left[\,\frac{\gamma_n}{\mu_n} + \left(\frac{\e_n}{\delta_n}\right)^2\frac{1}{\gamma_n\mu_n}\,\right]},
\end{align}
where the terms in the square parenthesis correspond to the energy of the pink and purple region, the green region, and the blue region respectively.

In the case $u_0(x,y)\in\{a(y), b(y)\}$ almost everywhere, we have that $\mu_n$ must be of order $1$.
Moreover, it is possible to see that for each choice of vanishing sequences $\{\eta_n\}_n,\{\gamma_n\}_n$, the first term is always of higher order than the last two. In particular, there is no rescaling that allows to \emph{see} that energy contribution (the bulk energy).
On the other hand, it is possible to make the last two terms of the same order if and only if $\eta_n=\gamma_n=\frac{\e_n}{\delta_n}$. Therefore, we guess that $\omega_n=\frac{\e_n}{\delta_n}$.

We now show that this scaling identifies functions $u_0\in\mathcal{R}$.
Dividing \eqref{eq:scaling} by $\frac{\e_n}{\delta_n}$ we get
\[
\purp{\frac{\e_n}{\delta_n}} + \green{\frac{1}{\mu_n}\left(\, \eta_n\frac{\delta_n}{\e_n}
    + \frac{1}{\eta_n}\frac{\e_n}{\delta_n} \right)}
+ \blue{\frac{1}{\mu_n}\left(\, \gamma_n\frac{\delta_n}{\e_n} + \frac{1} {\gamma_n}\frac{\e_n}{\delta_n} \right)},
\]
which is finite as $\e_n\to0$ if and only if $\mu_n$ is of order $1$ and $\gamma_n=\eta_n=\frac{\e_n}{\delta_n}$.\\

This is the heuristic reason to choose the scaling $\omega_n=\frac{\e_n}{\delta_n}$: it is expected to give equipartition of the surface energies and to select minimizers of $u_0$ of the form \eqref{eq:good_minimizers}.
We will rigorously prove in the next sections that indeed, this is the correct scaling.

%%%%%%%%%%%%%%%%%%%%%%%%%%%%%%%%%%%%%%%%%%%%
%%%%%%%%%%%%%%%%%%%%%%%%%%%%%%%%%%%%%%%%%%%%
%%%%%%%%%%%%%%%%%%%%%%%%%%%%%%%%%%%%%%%%%%%%
%%%%%%%%%%%%%%%%%%%%%%%%%%%%%%%%%%%%%%%%%%%%

\subsection{The limiting functional}

Motivated by the heuristics of the previous section, we introduce the new family of functionals.

\begin{definition}
For $n\in\N$, we define the functional $G_n: L^1(\Omega;\R^M)\to[0,+\infty]$ as
\[
G^1_n(u)\coloneqq \frac{\delta_n}{\e_n}G_n(u)
    =\int_{\Omega} \left[\, \frac{\delta_n}{\e_n}
		W\left( \frac{x}{\delta_n}, u(x) \right) 
		\,+ \e_n\delta_n|\nabla u(x)|^2 \,\right] \,dx.
\]
\end{definition}

\begin{remark}\label{rem:rewriting_Gn}
In the following, it is convenient to write the functional $G_n$ by using the undoflding operator.
For $n\in\N$, define the functional $\widetilde{G}_n: L^1(Q;\R^M)\to[0,\infty)$ as
\[
\widetilde{G}^1_n(v) \coloneqq \int_Q \left[\,
    \frac{\delta_n}{\e_n} W(y, v(y))
        + \frac{\e_n}{\delta_n} |\nabla v(y)|^2 
	\,\right] \,dy,
\]
and note that it is possible to write
\[
G^1_n(u)=\int_\o\, \widetilde{G}^1_n(\unf u(x,\cdot)) \,dx + R_n(u),
\]
where, recalling \eqref{eq:unfolding-operator}, we set
\begin{align*}
R_n(u) &\coloneqq \int_{\Lambda_{\delta_n}} \left[\, \frac{\delta_n}{\e_n}
		W\left( \frac{x}{\delta_n}, u(x) \right) 
		\,+ \e_n\delta_n |\nabla u(x)|^2 \,\right] \,dx \\
&\hspace{2cm} -|\Lambda_{\delta_n}|\e_n\delta_n \int_Q|\nabla a(y)|^2\,dy.
\end{align*}
In particular, we can see $G^1_n$ as a functional defined on $L^1(\Omega; L^1(Q;\R^M))$ as follows
\[
G^1_n(v) = \int_\o\, \widetilde{G}^1_n( v(x,\cdot)) \,dx
    + \widetilde{R}_n(v),
\]
where
\begin{align*}
\widetilde{R}_n(v) &\coloneqq \int_{\Lambda_{\delta_n}} \int_Q \left[\,
   \frac{\delta_n}{\e_n} W\left( y, v(x,y) \right) 
	\,+ \frac{\e_n}{\delta_n}|\nabla_y v(x,y)|^2 \,\right] \,dy\, dx \\
&\hspace{2cm} -\int_{\Lambda_{\delta_n}} \int_Q \frac{\e_n}{\delta_n} |\nabla a(y)|^2 \,dy\, dx.
\end{align*}
These representations will be useful in the rest of the paper.
\end{remark}

We now introduce the interfacial energy density of the limiting functional.
It is convenient to introduce the function $\chi:\R^M\to\{1,\dots,k\}$ defined as $\chi(y)\coloneqq i$ if $y\in E_i$. Note that $\chi\in SBV(\R^M)$ and its jump set $J_\chi$ corresponds to points $y\in \R^M$ such that there exist only two indexes $i\neq j$ with $y\in\partial E_i\cap\partial E_j$. It can be identified as $S_\chi\setminus J_\chi$, where $S_\chi$ is the set of singular points of $\chi$.

\begin{definition}\label{def:dw}
For $p,q,z_0\in\R^M$, consider the class
\begin{equation*}
\mathcal{A}(p,q,z_0) \coloneqq \left\{ \gamma \in W^{1,1}((-1,1); \R^M) : \gamma(-1) = p, \gamma(0)=z_0, \gamma(1) = q \right\}.
\end{equation*}
Define $\dw: \left[\, J_\chi\cup \left( \overline{Q}\setminus S_\chi \right) \,\right] \times\R^M\times\R^M \to[0, \infty)$ as
\begin{align*}
\dw(y,p,q)&\coloneqq \inf \Biggl\{\int_{-1}^0 2\sqrt{W_i(y, \gamma(t))}|\gamma'(t)|dt + \int_0^1 2\sqrt{W_j(y, \gamma(t))}|\gamma'(t)|dt \Biggr\}
\end{align*}
if $\chi^-(y)=i$ and $\chi^+(y)=j$, where the infimum is taken over points $z_0\in\R^M$, and over curves $\gamma \in \mathcal{A}(p,q,z_0)$.
\end{definition}

\begin{remark}
Note that in the case $\chi^-(y)=\chi^+(y)$, we have that
\[
\dw(y,p,q) = \inf \Biggl\{\int_{-1}^12\sqrt{W_i(y, \gamma(t))}|\gamma'(t)|\,dt \Biggr\},
\]
where the infimum is taken over curves $\gamma \in W^{1,1}((-1,1); \R^M)$ such that $\gamma(-1) = p$, and $\gamma(1) = q$.

In the case $\chi^-(y)\neq \chi^+(y)$ we cannot impose a priori the point $v_0\in\R^M$ where the profile will pass through at time $t=0$, and therefore we need to infimize also over that parameter.
\end{remark}

\begin{definition}\label{def:R}
Let
\begin{align*}
\widetilde{\mathcal{R}} &\coloneqq \Bigl\{ v\in L^1(\R^N;\R^M) : v \text{ is $Q$-periodic}, 
    v(y)\in\{a(y), b(y)\} \text{ a.e.},\\
&\hspace{8cm} v_{|Q}\in \bvloc(Q_0;\R^M) \Bigr\},
\end{align*}
where $Q_0\coloneqq Q\setminus\{x\in Q : a(x)=b(x)\}$, and define the class
\[
\mathcal{R} \coloneqq \left\{\, v\in L^1(\Omega;L^1(Q;\R^M)) \,:\, 
    \widetilde{v}(x,\cdot)\in \widetilde{\mathcal{R}} \text{ for a.e. } x\in\o \,\right\},
\]
where $\widetilde{v}:\R^N\to\R^M$ denotes the $Q$-periodic extension of $v\in L^1(Q;\R^M)$.
\end{definition}

We are now in position to define the limiting functional.

\begin{definition}
Let $G^1:L^1(\Omega;L^1(Q;\R^M))\to[0,+\infty]$ be defined as
\[
G^1(u)\coloneqq 
\left\{
\begin{array}{ll}
\displaystyle\int_\Omega \widetilde{G}^1(\widetilde{u}(x,\cdot)) \, dx & \text{ if } u\in \mathcal{R} ,\\
& \\
+\infty & \text{ else},
\end{array}
\right.
\]
where, for a function $v\in L^1(\R^N;\R^M)$, we set
\[
\widetilde{G}^1(v) \coloneqq \int_{\widetilde{Q}\cap J_v} \dw(y,v^-(y),v^+(y)) \dhno(y),
\]
and $\widetilde{Q}\coloneqq [0,1)^N$.
\end{definition}

\begin{remark}
Note that the energy $G^1$ is well defined. Indeed, for $u\in \mathcal{R}$, by using the measurability of $x\mapsto \widetilde{u}(x,\cdot)$, it is easy to see that the function $x\mapsto \widetilde{G}^1(\widetilde{u}(x,\cdot))$ is measurable.
Moreover, the jump set $J_v$ of a measurable function $v:\R^M\to\R$ is $\mathcal{H}^{M-1}$-rectifiable (see \cite{DelNin}).
Finally, as noted in \cite[Remark 1.8]{CriGra}, there are functions $v\in L^1(\R^N;\R^M)$ for which $\widetilde{G}^1(v)<\infty$, but $v\not\in \mathrm{BV}(\R^N;\R^M)$.
\end{remark}

The main result of this section is the following.

\begin{theorem}\label{thm:1_order}
Let $\{\e_n\}_{n}, \{\delta_n\}_{n}\subset (0,1)$ be infinitesimal sequences such that
\[
\lim_{n\to\infty}\frac{\e_n}{\delta_n} = 0.
\]
If $\{u_n\}_{n}\subset W^{1,2}(\o;\R^M)$ is such that
\[
\sup_n G^1_n(u_n) < +\infty,
\]
then there exists $u\in \mathcal{R}$ with $G^1(u)<\infty$ such that, up to a subsequence (not relabeled), $u_n\sts u$  strongly in $L^1(\Omega; L^1(Q;\R^M))$.
Moreover,
\[
G^1_n \stackrel{\Gamma}{\rightarrow} G^1
\]
with respect to strong two-scale convergence in $L^1(\Omega;L^1(Q;\R^M))$.
\end{theorem}

The result of Theorem \ref{thm:1_order} is written in the language of two-scale convergence. It gives a natural way to write the $\Gamma$-convergence result with respect to the strong $L^2$ covnergence.

\begin{corollary}\label{cor:main_result_x}
Let $\{\e_n\}_{n}, \{\delta_n\}_{n}\subset (0,1)$ be infinitesimal sequences such that
\[
\lim_{n\to\infty}\frac{\e_n}{\delta_n} = 0.
\]
Then
\[
G^1_n \stackrel{\Gamma}{\rightarrow} H^1
\]
with respect to the weak $L^2(\o;\R^M)$ convergence. Here
\[
H^1(v) \coloneqq \int_\Omega \widetilde{H}^1(v(x))\, dx,
\]
where
\[
\widetilde{H}^1(p) \coloneqq \min \left\{\, \widetilde{G}^1(w) \,:\, w\in \mathcal{R},\, \int_Q w(y)dy = p  \,\right\}.
\]
Moreover, if $\{u_n\}_{n}\subset W^{1,2}(\o;\R^M)$ is such that
\[
\sup_n G^1_n(u_n) < +\infty,
\]
then there exists $u\in \mathcal{R}$ with $H^1(u)<\infty$ such that, up to a subsequence (not relabeled), $u_n\rightharpoonup u$ weakly in $L^2(\Omega; \R^M)$.
\end{corollary}

Moreover, as for the case of the zeroth order $\Gamma$-limit, the mass constrain passes to the limit. Namely, the following holds.

\begin{corollary}\label{cor:main_result_mass}
Let $\{\e_n\}_{n}, \{\delta_n\}_{n}\subset (0,1)$ be infinitesimal sequences such that
\[
\lim_{n\to\infty}\frac{\e_n}{\delta_n} = 0.
\]
Fix $\overline{m}\in\R^M$, and define
\[
\mathcal{G}^1_n(u) \coloneqq
\begin{dcases} 
G^1_n(u;\overline{m}) & \text{ if } \int_\o u\, dx = \overline{m}, \\[5pt]
+\infty & \text{ else.} 
\end{dcases}
\]
Then it holds that $\mathcal{G}^1_n \stackrel{\Gamma}{\rightarrow} \mathcal{G}^1$ with respect to the strong two-scale convergence in $L^1(\Omega;L^1(Q;\R^M))$, where
\[
\mathcal{G}^1(u;\overline{m}) \coloneqq 
\begin{dcases} 
G^1(u;\overline{m}) & \text{ if } \int_\o\int_Q u\, dy\, dx = \overline{m}, \\[5pt]
+\infty & \text{ else.} 
\end{dcases}
\]
\end{corollary}

A similar result holds, with the obvious modifications, for the functionals considered in Corollary \ref{cor:main_result_x}.

Finally, we study the minimization problem for the limiting functional, with and without mass constrain. The proofs follow easily from the definition of the functionals, and by using a measurable selection result like that used in step 1 of the proof of Theorem \ref{thm:0_order}.

\begin{corollary}\label{cor:min_energy_first}
It holds that
\[
\min\left\{\, G^1(u) \,:\, u\in \mathcal{R} \,\right\}=0
\]
if and only if the $Q$-periodic extension of the whole $\R^N$ of the functions $a$ and $b$ are continuous.
Fix $\overline{m}\in\R^M$. Then
\[
\min\left\{\, \mathcal{G}^1(u;\overline{m}) \,:\, u\in \mathcal{R} \,\right\}=0
\]
if and only if the $Q$-periodic extension of the whole $\R^N$ of the functions $a$ and $b$ are continuous, and there exists $u\in\mathcal{R}$ with $\int_\o\int_Q u\, dy\,dx = \overline{m}$.
\end{corollary}

%%%%%%%%%%%%%%%%%%%%%%%%%%%%%%%%%%%%%%%%%%%%
%%%%%%%%%%%%%%%%%%%%%%%%%%%%%%%%%%%%%%%%%%%%
%%%%%%%%%%%%%%%%%%%%%%%%%%%%%%%%%%%%%%%%%%%%
%%%%%%%%%%%%%%%%%%%%%%%%%%%%%%%%%%%%%%%%%%%%

\subsection{Compactness}

This section is devoted to the proof of compactness, that we state separately.

\begin{lemma}\label{lem:compactness_special_class_0}
Let $\{u_n\}_n\subset W^{1,2}(\o;\R^M)$ be such that
\begin{equation}\label{eq:energy_bound_Gn}
\sup_n G^1_n(u_n) < +\infty.
\end{equation}
Then, up to a subsequence (not relabeled), $u_n\to u$ strongly two scale in $L^1(\Omega; L^1(Q;\R^M))$ for some $u\in \mathcal{R}$ with $G^1(u)<\infty$.
\end{lemma}

\begin{proof}For the sake of notation, we will write $\hat{u}_n$ in place of $\unf u_n$.\\

\textbf{Step 1.} Recalling Remark \ref{rem:rewriting_Gn}, we can write
\[
G^1_n(u_n) = \int_\o \widetilde{G}_n(\hat{u}_n(x,\cdot))\, dx + R_n(u_n),
\]
where
\begin{align*}
R_n(u_n) &\coloneqq \int_{\Lambda_{\delta_n}} \left[\, \frac{\delta_n^2}{\e_n}
		W\left( \frac{x}{\delta_n}, u_n(x) \right) 
		\,+ \e_n\delta_n|\nabla u_n(x)|^2 \,\right] \,dx \\
&\hspace{3cm}-|\Lambda_{\delta_n}|\e_n\delta_n \int_Q|\nabla a(y)|^2\,dy.
\end{align*}
Since
\[
R_n(u_n) \geq -|\Lambda_{\delta_n}|\e_n\delta_n \int_Q|\nabla a(y)|^2\,dy,
\]
and the right-hand side tends to zero as $n\to\infty$, from \eqref{eq:energy_bound_Gn} we get that
\[
\sup_n \int_\o \widetilde{G}_n(\hat{u}_n(x,\cdot))\, dx \leq C,
\]
for some $C<\infty$. \\

\textbf{Step 2.}
We claim that if is possible to find a subsequence $\{u_{n_j}\}_{j\in\N}$ such that, for $\mathcal{L}^N$-a.e. $x\in\o$, it holds
\begin{equation}\label{eq:bound_inside}
\limsup_{j\to\infty} \int_Q \left[\, \frac{\delta_n}{\e_n}
		W(y, \hat{u}_{n_j}(x,y))
		+ \frac{\e_n}{\delta_n} |\nabla_y \hat{u}_{n_j}(x,y)|^2 
	\,\right] \,dy < +\infty.
\end{equation}
For each $n\in\N$ define the function $f_n:L^1(\o)\to[0,+\infty]$ by
\[
f_n(x)\coloneqq \int_Q \left[\, \frac{\delta_n}{\e_n}
		W(y, \hat{u}_n(x,y))
		+ \frac{\e_n}{\delta_n} |\nabla_y \hat{u}_n(x,y)|^2 
	\,\right] \,dy.
\]
Then, by assumption, since $W\geq0$, we have $\sup_n\|f_n\|_{L^1(\o)}<+\infty$.
By the Chacon biting lemma (see \cite[Lemma 2.63]{FonLeo}), we have the following. 
There exists a subsequence $\{f_{n_k}\}_{k\in\N}$, and a sequence $\{r_{n_k}\}_{k\in\N}\subset(0,+\infty)$ with $\lim_{k\to\infty} r_{n_k}=+\infty$ such that, setting
\[
F_j \coloneqq \bigcup_{k=j}^\infty
    \left\{\, x\in\Omega \,:\, f_{n_k}(x)\geq r_{n_k}  \,\right\},
\]
we have $|F_j|\to0$ as $j\to\infty$.
Set
\[
 F\coloneqq \left\{ x\in\o \,:\, \limsup_{k\to\infty} f_{n_k}(x)=+\infty \right\}.
\]
Since it is possible to write
\[
F=\bigcap_{j\in\N} F_j,
\]
and $\{F_j\}_{j\in\N}$ is a decreasing sequence of sets contained in $\o$, we obtain that $|F|=0$.\\

\textbf{Step 3.} 
Let $x\in\o\setminus F$. Considering a sequence of compact sets $\{K_i\}_i$ invading $Q\setminus Q_0$, using \eqref{eq:bound_inside}, and \cite[Proposition 4.1]{CriGra} (see also the proof of Theorem 1.9 in \cite{CriGra}) we can extract a subsequence $\{u_{n_k}(x,\cdot)\}_{k\in\N}$ (possibly depending on $x\in\o\setminus F$), and find a function $v_x\in L^1(Q;\R^M)$ such that
\begin{itemize}
\item[(i)] $v_x(y)\in\{a(y), b(y)\}$ for a.e. $y\in Q$;
%\item[(ii)] the set $\{y\in Q : v_x(y)=a(y)\}$ has finite perimeter in $Q$;
\item[(ii)] $\hat{u}_{n_k}(x,\cdot)\to v_x$ strongly in $L^1(Q;\R^M)$ as $k\to\infty$;
\item[(iii)] $v_x\in BV_{loc}(Q\setminus Q_0;\R^M)$;
\item[(iv)] $\widetilde{G}^1(v_x)<\infty$.
\end{itemize}
We want to prove that the subsequence does not depend on the point $x\in\o\setminus F$.
Note that \eqref{eq:energy_bound_Gn} implies that
\[
\sup_n G_n(u_n)<\infty,
\]
and thus Theorem \ref{thm:0_order} gives the existence of a subsequence $\{u_{n_j}\}_{j}$ and of a function $\widetilde{u}\in L^2(\o;\R^M)$ such that $u_{n_j}\rightharpoonup \widetilde{u}$ in $L^2(\o;\R^M)$.
In particular, since
\begin{equation}\label{eq:bond_uj_L2}
\sup_{j\in\N} \| u_{n_j} \|_{L^2(\o;\R^M)} < \infty,
\end{equation}
by applying Proposition \ref{prop:2-scale-compactness}, we get that there exists a (not relabeled) subsequence such that $u_{n_j}\rightharpoonup u$ weakly two-scale in $L^2(\o;L^2(\o;\R^M))$, for some $u\in L^2(\o;L^2(\o;\R^M))$. Therefore
\[
\widetilde{u}(x) = \int_Q u(x,y)\, dy,
\]
for a.e. $x\in\o$, and, by using (iii), it easy to see that $u_{n_j}(x,\cdot)\to u(x,\cdot)$ strongly in $L^1(Q;\R^M)$ for all $x\in\o\setminus F$.

Finally, we claim that $u_{n_j}\to u$ strongly two scale in $L^1(\o;L^2(Q;\R^M))$. 
Define, for each $j\in\N$, $g_j:\o\to [0,\infty)$ as
\[
g_j(x) := \|\hat{u}_{n_j}(x,\cdot) - u(x,\cdot)\|_{L^1(Q;\R^M)}.
\]
Then from \eqref{eq:bond_uj_L2} we get
\[
\sup_{j\in\N}\int_\o {g_j^2\;dx}<+\infty.
\]
Using De la Val\'ee Poussin criteria, we have that $\{g_j\}_{j}$ is equiintegrabile. Furthermore, by (iii), we get that $g_j \to 0$ pointwise almost everywhere.
We can now apply Vitali Convergence Theorem to conclude that $g_j\to 0$ in $L^p(\o)$ strong for any $p\in[1,2)$. This concludes the proof of the compactness result.
\end{proof}

%%%%%%%%%%%%%%%%%%%%%%%%%%%%%%%%%%%%%%%%%%%%
%%%%%%%%%%%%%%%%%%%%%%%%%%%%%%%%%%%%%%%%%%%%
%%%%%%%%%%%%%%%%%%%%%%%%%%%%%%%%%%%%%%%%%%%%
%%%%%%%%%%%%%%%%%%%%%%%%%%%%%%%%%%%%%%%%%%%%

\subsection{Liminf inequality}

The main result of this section is the following.

\begin{proposition}\label{prop:liminf}
Let $u\in L^1(\o; L^1(Q;\R^M))$ and let $\{u_n\}_{n}\subset W^{1,2}(\o;\R^M)$ with $u_n\to u$ strongly two scale in $L^1(\o; L^1(Q;\R^M))$. Then
\[
G^1(u) \leq \liminf_{n\to\infty} G^1_n(u_n).
\]
\end{proposition}

The proof of Proposition \ref{prop:liminf} is based on the liminf inequality for a single periodicity cell $Q$. This result is essentially contained in \cite[Proposition 4.2]{CriGra} and in the remarks made in the proof of \cite[Theorem 1.9]{CriGra}. In the language of this paper, it writes as follows.

\begin{proposition}\label{prop:liminfQ}
Let $v\in L^1(Q;\R^M)$ such that $v(y)\in\{a(y), b(y) \}$ for a.e. $y\in Q$. Let $\{v_n\}_n\subset W^{1,2}(Q;\R^M)$ with $v_n\to v$ in $L^1(Q;\R^M)$. Then
\[
\widetilde{G}^1(v) \leq \liminf_{n\to\infty} \widetilde{G}_n(v_n).
\]
\end{proposition}

The proof of Proposition \ref{prop:liminfQ} we present here uses a slightly different strategy from that of \cite[Proposition 4.2]{CriGra}, and requires some technical results. We decided to show how to get Proposition \ref{prop:liminf} once Proposition \ref{prop:liminfQ} is established, and then to move to the technical results needed to obtain this latter.

\begin{proof}[Proof of Proposition \ref{prop:liminf}]
Let $u\in L^1(\o;L^1(Q;\R^M))$, and take $\{u_n\}_n\subset W^{1,2}(\o;\R^M)$ such that $u_n\to u$ strongly two-scale in $L^1(\o;L^1(Q;\R^M))$.
Without loss of generality, we can assume that
\[
\liminf_{n\to\infty} G_n(u_n)<+\infty,
\]
otherwise there is nothing to prove. By the compactness result (see Lemma \ref{lem:compactness_special_class_0}), we get that $u\in\mathcal{R}$.
Therefore, recalling the arguments in the proof of Lemma \ref{lem:compactness_special_class_0} we get that
\begin{align*}
\liminf_{n\to\infty} G_n(u_n) &\geq
    \liminf_{n\to\infty} \int_\o \widetilde{G}_n(T_{\delta_n} u_n(x,\cdot))\, dx + R_n(u_n) \\
&= \liminf_{n\to\infty} \int_\o \widetilde{G}_n(T_{\delta_n} u_n(x,\cdot))\, dx \\
&\geq \int_\o \liminf_{n\to\infty} \widetilde{G}_n(T_{\delta_n} u_n(x,\cdot))\, dx \\
&\geq \int_\o G^1(u(x,\cdot))\, dx,
\end{align*}
where the previous to last step follows by the Fatou's lemma, while last step is justified by the fact that $\unf u_n(x,\cdot)\to u(x,\cdot)$ for a.e. $x\in\o$, together with Proposition \ref{prop:liminfQ}. This concludes the proof of the liminf inequality.
\end{proof}

%%%%%%%%%%%%%%%%%%%%%%%%%%%%%%%%%%%%%%%%%%%%
%%%%%%%%%%%%%%%%%%%%%%%%%%%%%%%%%%%%%%%%%%%%
%%%%%%%%%%%%%%%%%%%%%%%%%%%%%%%%%%%%%%%%%%%%
%%%%%%%%%%%%%%%%%%%%%%%%%%%%%%%%%%%%%%%%%%%%

\subsubsection{Bound on the Euclidean length of geodesics}\label{sec:Euclidean}

First, we prove a technical lemma on bounds of Euclidean length of geodesics necessary for the liminf inequality in $Q$. While the overall proof strategy is similar to that in \cite{CriGra} and \cite{sternberg88}, our construction by estimating the energy within each level set (see Step 2 of \Cref{lem:bded_eu_len_well}) is novel.

In this section, in order to make the notation lighter, we will make the following abuse of notation. Fix a vector $\nu\in\S^{N-1}$, a point $x_0\in Q$, and a unit square $C$ centered at the origin and with two faces orthogonal to $\nu$.
For $t>0$, we denote by $Q'_t\coloneqq (x_0 + t C)\cap \nu^\perp$.
For $y'\in Q'_t$ and $z\in\R$, we denote point $y'+z\nu$ by $(y',z)$.
Note that the fact that some of the above points could be outside $\Omega$ is of no concern for us, since all of the functions that we consider are $Q$ periodic, and thus can be naturally extended from $\Omega$ to the whole $\R^N$.

We are now in position to introduce the minimization problem that will be investigated in this section.

\begin{definition}
For $p,q \in \R^M$, let
\begin{equation*}
\mathcal{A}(p,q) \coloneqq \left\{ \gamma \in W^{1,1}([-1,1]; \R^M) : \gamma(-1) = p\;\;\textup{and}\;\;\gamma(1) = q \right\}.
\end{equation*}
For $\e>0$, $y'\in Q'_\e$, and $p, q\in\R^M$, define
\[
H_\ep(y',p,q)\coloneqq  \inf\left\{\int_{-1}^{1} {F_\ep(y',\gamma(t))|\gamma'(t)| dt}\; : \gamma \in \mathcal{A}(p,q)\right\},
\]
where
\[
F_\ep(y',p) \coloneqq  \underset{|z|\leq \ep}{\min}{\sqrt{W(y',z,p)}}.
\]
\end{definition}

The main result of this section is the following.

\begin{theorem}\label{thm:boundEL}
Fix $x_0\in Q$, $\nu\in\S^{N-1}$, and $\widetilde{R}>0$.
If $x_0\in\cup_{i=1}^k \partial E_i$, assume that it belongs to only one of those sets.
Then there exist $\e_1>0$ and $L>0$ such that, given any $\e\in(0,\e_1)$, $y'\in Q'_\e$, and $p,q\in B(0,\widetilde{R})$, the minimization problem defining $H_\ep(y',p,q)$ admits a solution $\gamma\in \mathcal{A}(p,q)$ such that
\[
\int_{-1}^1 |\gamma'(t)| dt \leq L.
\]
\end{theorem}

The strategy to prove Theorem \ref{thm:boundEL} is the following.
First we consider the case where the point $x_0\in Q\setminus \cup_{i=1}^k \partial E_i$. This means that $W=W_i$ for some $i\in\{1,\dots,k\}$, and, in particular, it is Lipschitz in the second variable.
Consider a minimizing sequence $\{\gamma_j\}_{j\in\N}\subset\mathcal{A}(p,q)$ for the minimization problem defining $H_\e(y',p,q)$. The first step is to investigate the behaviour of the sequence of curves, for $\e$ sufficiently small, close by and far away from the wells.
In particular, in Lemma \ref{lem:bded_eu_len_well} we prove that the portion of the curves $\gamma_j$ that is sufficiently close to $a(y',0)$ (or to $b(y',0)$) has uniformly bounded Euclidean length.
Then, by using a lower bound on $F_\e$ far from the wells, we conclude that also the Euclidean length of the sequence $\{\gamma_j\}_{j\in\N}$ is uniformly bounded in that region.
The proof of Theorem \ref{thm:boundEL} then follows by using a standard argument based on the Ascoli-Arzel\`{a} Theorem ensuring the existence of minimizing geodesics for the minimization problem defining $H_\e(y',p,q)$. We refer to (see \cite[Lemma 3.1]{CriGra} for details.
Finally, the case $x_0\in \cup_{i=1}^k \partial E_i$ may be deduced from the previous case.\\

We start by collecting some basic properties of $F_\e$, which can verified easily from the definition. In particular, for a proof of (3), we refer to \cite[Proposition 3.2]{CriGra}.

\begin{proposition}[Properties of $F_\ep$]\label{prop:Prop_F}
The followings hold:
\begin{enumerate}
    \item The function $(y',p)\mapsto F_\ep$ is Lipschitz;
    \item We have $F_\ep(y',p) = 0$ if and only if $p \in \{a(y',t),b(y',t)\}$ for some $t \in [-\ep,\ep]$;
    \item If there is a function $g: S \to \R^M$ such that
    \[
    W(y',y_n,p) = |p-g(y',y_n)|^2,
    \]
    then $F_\ep(y',p) = \dist(p,\gr_\e^g(y'))$, where
    \[
    \gr_\e^g(y')\coloneqq  \{g(y',t) \;:\; y'\in Q',  \;|t|\leq \e \}.
    \]
\end{enumerate}
\end{proposition}
 
Next, we state a property, based on the parametrization invariant characteristic of the minimization problem defining $H_\varepsilon$, that will be used several times.

\begin{proposition}\label{prop:Prop_H}
Let $p,q \in \R^M$, $y' \in Q'_\e$,
and $\gamma \in \mathcal{A}(p,q)$ such that 
\[
\int_{-1}^{1} {F_\ep(y',\gamma (t))|\gamma'(t)| dt}
    \leq H_\ep(y',p,q)+\frac{1}{j},
\]
for some $j\in\N$.
Then
\[
\int_{t_1}^{t_2} {F_\ep(y',\gamma(t))|\gamma'(t)| dt}
    \leq H_\ep(y',\gamma(t_1),\gamma(t_2))+\frac{1}{j},
\]
for all $[t_1,t_2] \subset [-1,1]$.
\end{proposition}

The main idea in the proof of the bound of the Euclidean length close by the wells is to consider level sets of $F(y',\gamma(\cdot))$ in the construction of a competitor for the minimization problem defining $H_\e(y',p,q)$.

\begin{definition}
For $\e>0$, $y' \in Q'_\e$, $p,q\in\R^M$, and $\gamma\in\mathcal{A}(p,q)$, and $k\in\N\setminus\{0\}$, we define
\[
T^{k}_\e(y',\gamma) \coloneqq  \left\{t\in [-1,1]\;:\;\frac{1}{(k+1)^2}<F_\e(y',\gamma(t))\leq \frac{1}{k^2}\right\}.
\]
\end{definition}

\begin{remark}
Note that by continuity of $F_\e$, we have that
\[
\overline{T^{k}_\e(y',\gamma)} = \left\{t\in [-1,1]\;:\;\frac{1}{(k+1)^2}\leq F_\e(y',\gamma(t))\leq \frac{1}{k^2}\right\}.
\]
\end{remark}

Now we are ready to prove the key technical lemma of this section.

\begin{lemma}\label{lem:bded_eu_len_well}
Let $x_0\in Q\setminus \cup_{i=1}^k \partial E_i$ and $\nu\in\S^{N-1}$.
Then there exist $r_0>0$, $\e_0 > 0$, and $L_1>0$, such that for any $0<\e < \ep_0$, $y' \in Q'_\e$, $p,q \in B_{r_0}(y')\coloneqq B_{r_0}(a(y',0))$, the following property holds.
Let $\{\gamma_j\}_j \subset \mathcal{A}(p,q)$ be a minimizing sequence for the minimization problem defining $H_\e(y',p,q)$. Then
\[
\int_{T_j} |\gamma'_j(t)|dt \leq L_1,
\]
where $T_j\coloneqq \{t\in[-1,1] \,:\, \gamma_j(t)\in B_{r_0}(y') \}$. The same result holds with $B_{r_0}(y')\coloneqq B_{r_0}(b(y',0))$.
\end{lemma}

\begin{proof}
Let $R>0$ be as in \ref{W4}.
We define
\[
\e_0 \coloneqq \min\left\{\frac{R}{2\mathrm{Lip}(a)}, \mathrm{dist}(x_0,\cup_{i=1}^k \partial E_i) \,\right\},
\]
where $\mathrm{Lip}(a)$ denotes the maximum over the index $i$ of the Lipschitz constant of $a_i$ in $E_i$.
Moreover, by using the uniform lower bound on $W(x,\cdot)$ close to the wells, it is possible to choose $r_0>0$ so small such that
\[
B_{r_0}(y') \subset \bigcup_{k=1}^\infty \left\{\, p\in\R^M \,:\,
    \frac{1}{(k+1)^2} < F_\e(y',p) \leq \frac{1}{k^2}  \,\right\},
\]
for all $\e\in(0,\e_0)$ and all $y'\in Q'_\e$.
Up to further decreasing the value of $r_0$, we can also suppose that $r_0\leq R/2$.\\

\underline{Case 1.} Assume that $H_\ep(y',p,q)=0$. Note that this happens if and only if $p,q\in Z_\e \coloneqq \{ a(y',t) : t\in[-1,1] \}$. In this case, the solution to the minimization problem defining $H(y', p, q)$ is given by the curve lying in $Z_\e$ joining the two points. By using the Lipschitz regularity of $a$, we get that its Euclidean length is less than $2\mathrm{Lip}(a)\e_0$.\\

\underline{Case 2.} Assume that $H_\ep(y',p,q)>0$. Without loss of generality, we can assume that
\begin{equation}\label{eq:almost_min}
    \int_{-1}^{1} F_\ep(y',\gamma_j(t))|\gamma'_j(t)| dt \leq 2 H_\ep(y',p,q),
\end{equation}
for all $j\in\N$.

\textbf{Step 1: Bounds on $F_\e$.}
We claim that if $\e\in(0,\e_0)$ and $r\in(0,r_0)$, for all $y'\in Q'_\e$ and $z \in B_r(y')$ it holds
\begin{equation}\label{eq:quad_bounds}
    \frac{1}{\sqrt{c_1}}\dist(z,\gr_\e^a(y'))\leq F_\e(y',z)\leq\sqrt{c_1}\dist(z,\gr_\e^a(y')).
\end{equation}
Indeed, by the triangle inequality and Lipschitz regularity of $a$, we get for $p \in B_{r_0}(y')$:
\[
|p-a(y',\e t)|\leq |p - a(y',0)|+\textup{Lip(a)}\e |t|< r_0 + \textup{Lip(a)}\e_0 \leq  R.
\]
This inequality gives that $B_{r_0}(y') \subset T_{R}(\gr_\e^a(y'))$, and the desired inequality follows by applying Property 2 in Proposition \ref{prop:Prop_F}.

From \eqref{eq:quad_bounds}, we easily have the bound $F_\e(y',z)\leq \sqrt{c_1}R$ for $z \in B_{r_0}(y')$.\\

\begin{figure}
    \centering
    \includegraphics[scale=0.5]{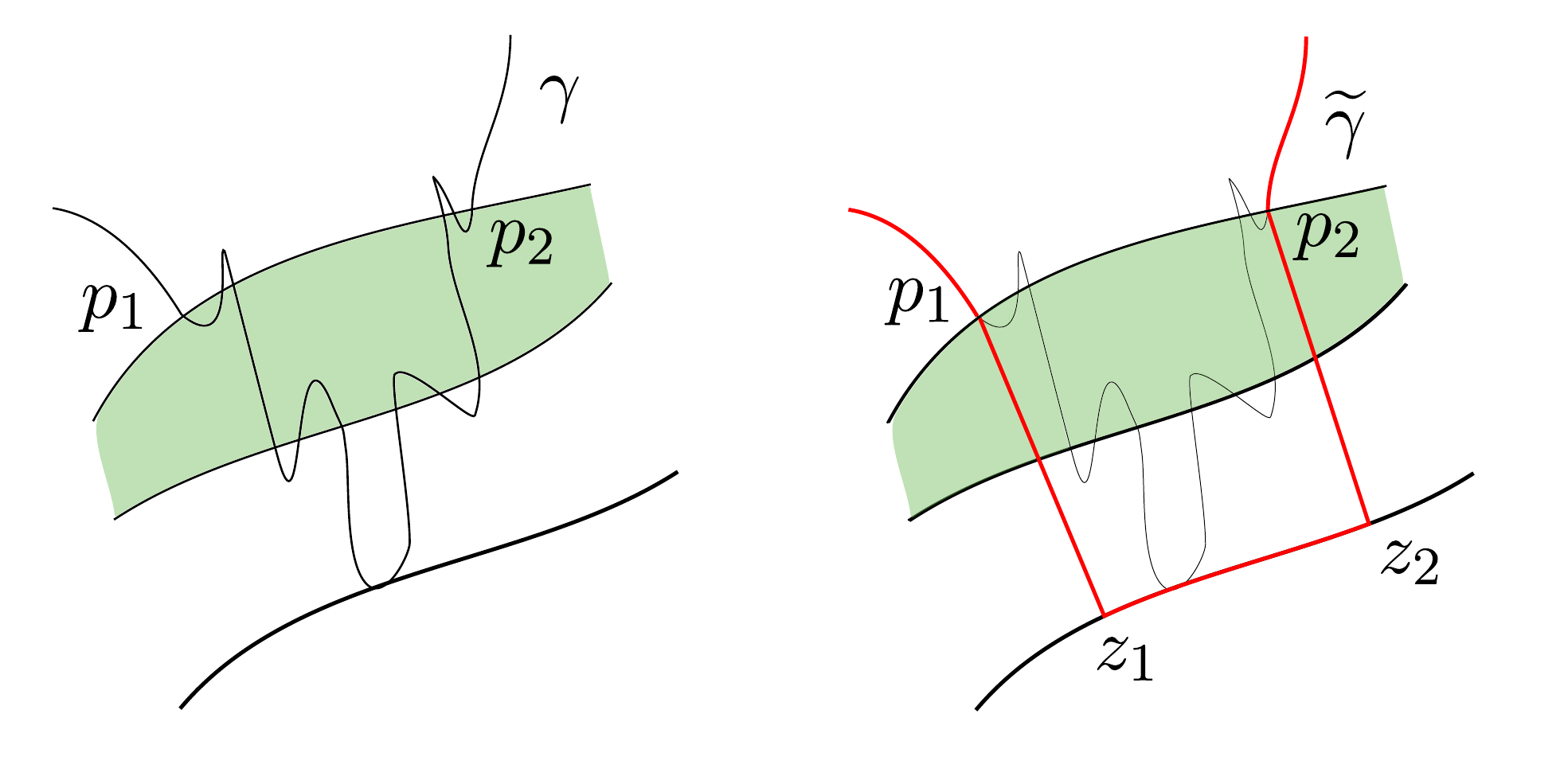}
    \caption{The construction of the competitor $\widetilde{\gamma}$ carried out in step 2, in order to estimate the length of a curve $\gamma_j$ in the green region, namely in the region between two level sets of $F_\varepsilon(y',\cdot)$.}
    \label{fig:competitor}
\end{figure}

\textbf{Step 2: Bounding the energy in a level set. }
Fix $j, k\in\N\setminus\{0\}$. Without loss of generality, we will suppose that
\[
\gamma_j\left(T^k_\e(y',\gamma_j)\right)\cap B_{r_0}(y')\neq \emptyset.
\]
Furthermore, to ease notation, we will set
\[
\mathcal{T}^k_\ep\coloneqq  \left\{ t\in T^k_\e(y',\gamma_j) \,:\, \gamma_j(t)\in B_{r_0}(y')\right\}.
\]
We want to get a uniform bound on the Euclidean length of the curves $\gamma_j$'s in the set $T_\ep(y')$. Let
\[
t_1 \coloneqq  \inf\left\{t:\;t \in \mathcal{T}^k_\ep\right\}, \quad\quad\quad
t_2 \coloneqq  \sup\left\{t :\; t \in  \mathcal{T}^k_\ep\right\},
\]
and
\[
p_1 \coloneqq  \gamma_j(t_1), \quad\quad\quad
p_2 \coloneqq  \gamma_j(t_2).
\]
Using the definition of $L_\e(y')$ and the choice of $t_1,t_2$, we get the following simple lower bound on the energy:
\begin{align}\label{eq:lower_bd}
    \frac{1}{(k+1)^2}\int_{\overline{\mathcal{T}^k_\ep}} |\gamma'_j(t)| dt
    &\leq\int_{\overline{\mathcal{T}^k_\ep}} F_\ep(y',\gamma_j(t))|\gamma'_j(t)| dt \nonumber\\
    &\leq \int_{t_1}^{t_2} F_\ep(y',\gamma_j(t))|\gamma'_j(t)| dt .
\end{align}
We employ the geometric property of $H_\e$ given in \Cref{prop:Prop_H} and \eqref{eq:quad_bounds} to deduce that for any $\gamma \in \mathcal{A}(p_1,p_2)$, we have the bound
\begin{align}\label{eq:truncated_geodesic}
    \int_{t_1}^{t_2} {F_\ep(y',\gamma_j(t))|\gamma'_j(t)| dt}
    &\leq 2H_\e(y',p_1,p_2) \nonumber\\
    &\leq 2\int_{-1}^{1} {F_\ep(y',\gamma(t))|\gamma'(t)| dt} \nonumber \\
    &\leq  2\int_{-1}^{1} {\sqrt{c_1}\dist(\gamma(t),\gr_\e^a(y'))|\gamma'(t)| dt}.
\end{align}
In order to further bound from above the right-hand side of the above expression, we will construct a suitable competitor $\tilde{\gamma}\in\mathcal{A}(p+1,p_2)$.
First, note that there exist $z_1, z_2 \in \gr_\e^a(y')$ such that
\[
\dist(p_1,\gr_\e^a(y')) = |p_1- z_1|,\quad\quad\quad
\dist(p_2,\gr_\e^a(y')) = |p_2- z_2|.
\]
We define $\tilde{\gamma}\in\mathcal{A}(p_1, p_2)$ as the union of the following three curves (see Figure \ref{fig:competitor}):
\begin{enumerate}
    \item The segment between $p_1$ and $z_1$;
    \item The portion of $\gr_\e^a(y')$ that connects the points $z_1$ and $z_2$;
    \item The segment between $z_2$, and $p_2$.
\end{enumerate}
Note that, since the energy is parameterization invariant, we do not have to specify the precise parametrization of $\tilde{\gamma}$.
We now estimate the energy of $\tilde{\gamma}$. Again, by the parametrization invariant property of the functional, we can use $\pm1$ as initial and final time respectively for each of the three curves. Note only the two segments contributes to the energy of the curve, and by a direct evaluation, we get that
\begin{align}\label{eq:upper_bd}
\int_{-1}^{1} {\dist(\tilde{\gamma}(t),\gr_\e^a(y'))|\tilde{\gamma}'| dt}
&= \frac{\dist^2(p_1,\gr_\e^a(y')) + \dist^2(p_2,\gr_\e^a(y')}{2} \nonumber \\
&\leq \max\left\{\, \dist^2(p_1,\gr_\e^a(y')), \dist^2(p_2,\gr_\e^a(y')  \,\right\} \nonumber \\ &\leq c_1 \max\left\{\, \left( F_\e(y',p_1) \right)^2, \left( F_\e(y',p_2) \right)^2  \right\} \nonumber \\
&\leq \frac{c_1^2}{k^4},
\end{align}
where the previous to last inequality follows from \eqref{eq:quad_bounds}, while last step is justified by the fact that, since $p_1, p_2 \in \gamma_j(\overline{T_\e(y')})$, it holds
\[
\max\{F_\e(y',p_1),F_\e(y',p_2)\} \leq \frac{1}{k^2}.
\]

Thus, combining \eqref{eq:lower_bd}, \eqref{eq:truncated_geodesic}, and \eqref{eq:upper_bd}, we get
\[
\frac{1}{(k+1)^2}\int_{\overline{\mathcal{T}^k_\ep}}|\gamma'_j(t)| dt
    \leq \frac{2c_1^{\frac{3}{2}}}{k^4},
\]
that yields the upper bound
\begin{equation}\label{eq:lim_sup_level_set_energy}
\int_{\overline{\mathcal{T}^k_\ep}}|\gamma'_j(t)| dt \leq \frac{2c_1^{\frac{3}{2}}(k+1)^2}{k^4},
\end{equation}
for all $j\in\N$.\\

\textbf{Step 3: Bounding the Euclidean length.}
We have that
\begin{align}
    \int_{T_j} |\gamma'_j(t)| dt &\leq \left[\, \int_{Z_\e}|\gamma'_j(t)| dt
    + \sum_{k=1}^\infty \int_{\mathcal{T}^k_\ep}|\gamma'_j(t)| dt \,\right] \nonumber \\
    &\leq 2\,\textup{Lip}(a)\e_0 + \sum_{k=1}^\infty \frac{2c_1^{\frac{3}{2}}(k+1)^2}{k^4} =: L_1,
\end{align}
where last step follows from \eqref{eq:lim_sup_level_set_energy} and the estimate obtained in case 1. Note that the right-hand side is independent of $\e, y'$, and thus we achieve the desired result. 
\end{proof}

Now we are ready to prove the main result of this section.

\begin{proof}[Proof of Theorem \ref{thm:boundEL}]
Fix $x_0\in Q$, $\nu\in\S^{N-1}$, $\widetilde{R}>0$. Take $p,q\in B(0,\widetilde{R})$.

\underline{Case 1} Assume $x_0\in Q\setminus \cup_{i=1}^k \partial E_i$.

\textbf{Step 1: Choice of $\e_1$.}
Let $r_0, \e_0$ be given by Lemma \ref{lem:bded_eu_len_well}. Define
\[
\e_1 \coloneqq  \min \left\{\e_0, \frac{r_0}{2 \mathrm{Lip}(a)},\frac{r_0}{2 \mathrm{Lip}(b)} \right\}.
\]
Note that for every $z \notin \overline{B_{r_1}(b(y',0))} \cup \overline{B_{r_1}(a(y',0))}$ and for every $|t|\leq \e<\e_1$ we have
\begin{align}\label{eq:far_away_1}
    |z-a(y',t)|\geq |z-a(y',0)|-\mathrm{Lip}(a)|t|> r_0-\mathrm{Lip}(a)\e_2\geq \frac{r_0}{2},\\
    |z-b(y',t)|\geq |z-b(y',0)|-\mathrm{Lip}(b)|t|> r_0-\mathrm{Lip}(b)\e_2\geq \frac{r_0}{2}.\label{eq:far_away_2}
\end{align}
Fix $0<\e<\e_1$ and $y'\in Q'_\e$.\\

\textbf{Step 2: Estimate of the Euclidean length.}
Let $\{\gamma_j\}_j\subset \mathcal{A}(p,q)$ be a sequence satisfying \eqref{eq:almost_min}.
We will bound the Euclidean length of each $\gamma_j$ in each of the following regions separately:
\[
\mathcal{R}^1_j \coloneqq \{ t\in[-1,1] \,:\, \gamma_j(t)\in \overline{B_{r_0}(a(y',0))}  \},
\]
\[
\mathcal{R}^2_j \coloneqq \{ t\in[-1,1] \,:\, \gamma_j(t)\in \overline{B_{r_0}(b(y',0))}  \},
\]
and
\[
\mathcal{R}^3_j \coloneqq \{ t\in[-1,1] \,:\, \gamma_j(t)\in\overline{B_{r_0}(a(y',0))}^c \cap \overline{B_{r_0}(b(y',0))}^c \}.
\]
We start from the latter region. By \eqref{eq:far_away_1} and \eqref{eq:far_away_2}, together with Remark \ref{rem:bounds}, we get that there is $C_\frac{r_0}{2}>0$ such that
\[
F_\e \geq \sqrt{C_\frac{r_0}{2}}.
\]
Therefore
\begin{align}\label{eq:est_R3}
    \int_{\mathcal{R}^3_j} |\gamma' _j(t)|
    &\leq \frac{1}{\sqrt{C_\frac{r_0}{2}}}\int_{R^3} F_\e(y',\gamma_j(t))|\gamma' _j(t)| \nonumber \\
    &\leq \frac{1}{\sqrt{C_\frac{r_0}{2}}} \left(H_\ep(y',p,q)+\frac{1}{j} \right) \nonumber \\
    &\leq\frac{1}{\sqrt{C_\frac{r_0}{2}}} \left(m|p-q|+\frac{1}{j} \right) \nonumber \\
    &\leq \frac{1}{\sqrt{C_\frac{r_0}{2}}}(2m\widetilde{R}+1),
\end{align}
where
\[
m \coloneqq \sup\left\{W(y',z,p): y' \in Q'_{\e_2},\, z \in [-\e_2,\e_2],\, p\in B(0,\widetilde{R}) \right\} < \infty,
\]
and the previous to last step follows by considering as a competitor the segment joining $p$ and $q$.

We now bound the Euclidean length in the regions $\mathcal{R}^1_j$ and $\mathcal{R}^2_j$.
By using Proposition \ref{prop:Prop_H}, together with Lemma \ref{lem:bded_eu_len_well}, we obtain that
\begin{equation}\label{eq:est_R12}
\int_{\mathcal{R}^1_j} |\gamma'_j(t)| dt + \int_{\mathcal{R}^2_j} |\gamma'_j(t)| dt \leq 2L_1,
\end{equation}
where the constant $L_1<\infty$ depends only on $x_0$ and $\nu$.\\

\textbf{Step 3: Existence of a geodesic.} Using \eqref{eq:est_R3} and \eqref{eq:est_R12}, there exists $L>0$ depending only on $\widetilde{R}, x_0$, and $\nu$, such that
\[
\int_{-1}^1 |\gamma'_j(t)|\, dt\leq L
\]
for all $j\in\N$. A standard argument based on the Ascoli-Arzel\`a Theorem (see \cite[Lemma 3.1]{CriGra} for more details) yields the desired result.\\

\underline{Case 2.} 
Assume that $x_0\in \partial E_i\setminus \cup_{j\neq i}\partial E_j$, and that $\chi^-(y)=i$, $\chi^+(y)=j$.
It is easy to see that there exist $S>0$, depending only on $\widetilde{R}, x_0$, and $\nu$ such that, for every $0<\e<\e_1$ and $y'\in Q'_{\e}$ it is possible to find $z_0\in B(0,S)$ and curves $\gamma_i, \gamma_j\in W^{1,1}([-1,1];\R^M)$, with
\[
\gamma_i(-1)=p,\quad\quad \gamma_i(1)=z_0,\quad\quad
\gamma_j(-1)=z_0,\quad\quad \gamma_j(1)=q,
\]
such that (recall that the functional is invariant by reparametrization)
\[
H_\e(y',p,q) = \int_{-1}^1 F^i_e(y', \gamma_i(t)) |\gamma'_i(t)|\, dt + 
    \int_{-1}^1 F^j_e(y', \gamma_j(t)) |\gamma'_i(t)|\, dt,
\]
where
\[
 F^i_\e(y', s) = \min_{z\in[-\varepsilon,0]} \sqrt{W_i(y',z,s)},\quad\quad
 F^j_\e(y', s) = \min_{z\in[0,\varepsilon]} \sqrt{W_j(y',z,s)}.
\]
Thus, by applying case 1 to $\gamma_i$ and $\gamma_j$, and $S>0$, we conclude also in this case.
\end{proof}

We are now in position to prove the liminf inequality in $Q$. Since the strategy follows a similar argument to that of \cite[Proposition 4.2]{CriGra}, we will sketch the main ingredients of the proof, focusing on the points where the two arguments differ.

\begin{proof}[Proof of Proposition \ref{prop:liminfQ}]
Let $\{v_n\}_n\subset W^{1,2}(Q;\R^M)$ with $v_n\to v$ in $L^1(Q)$. Without loss of generality, we can assume that
\[
\limsup_{n\to\infty} \widetilde{G}_n(v_n)<\infty.
\]
Note that, thanks to assumption \ref{W4}, we can use the compactness argument in \cite[Proposition 4.1, and Theorem 1.9]{CriGra} to get that $v\in \widetilde{\mathcal{R}}$ with $\widetilde{G}^1(v)<\infty$.
Fix $x_0\in J_v\cap K$, where $K\subset Q_0$ is a compact set. The idea is to use a blow-up argument, as in the proof of \cite[Proposition 4.2]{CriGra}.
Note that, thanks to the continuity of the wells $a_i$'s and $b_i$'s, together with the fact that the blow up is a local argument, and thanks to assumption (W4), we can use of the estimate on the Euclidean length of solutions to the minimization problem defining $\dw$ provided by Theorem \ref{thm:boundEL}.
The only difference with the argument used in the proof of \cite[Proposition 4.2]{CriGra} is in step 2, where the functional $F_m$ (see equation (100) in \cite{CriGra}) is defined here as
\[
F_m(x', p)\coloneqq \inf\left\{\, \sqrt{W( (x',g_m(t)), p)} \,:\, t\in(-1,1)  \,\right\},
\]
where $g_m(t)\coloneqq (y',z+t\nu)$, for some $\nu\in\S^{N-1}$. Using Theorem \ref{thm:boundEL}, we obtain that the Euclidean length of the solutions to the geodesic problem
\[
\min\left\{\, \int_{-1}^1 F_m(x',\gamma(t)) |\gamma'(t)|\, dt  \,\right\},
\]
where the infimum is taken over $\gamma\in W^{1,1}((-1,1);\R^M)$ with $\gamma(-1)=p$ and $\gamma(1)=q$, are uniformly bounded with respect to $\varepsilon$ and $x'$.
We now have all the elements that allows us to conclude by following the same strategy.
\end{proof}

%%%%%%%%%%%%%%%%%%%%%%%%%%%%%%%%%%%%%%%%%%%%
%%%%%%%%%%%%%%%%%%%%%%%%%%%%%%%%%%%%%%%%%%%%
%%%%%%%%%%%%%%%%%%%%%%%%%%%%%%%%%%%%%%%%%%%%
%%%%%%%%%%%%%%%%%%%%%%%%%%%%%%%%%%%%%%%%%%%%

\subsection{Limsup inequality}

This section is devoted to the construction of the recovery sequence.

\begin{proposition}\label{prop:limsup}
Let $u\in L^1(\o; L^1(Q;\R^M))$. Then there exists a sequence $\{u_n\}_n\subset W^{1,2}(\o;\R^M)$ with $u_n\to u$ strongly two scale in $L^1(\o; L^1(Q;\R^M))$ such that $G_n(u_n)\to G^1(u)$ as $n\to\infty$.
\end{proposition}

The construction of the recovery sequence will be done in three steps: first for the class of simple functions in $\mathcal{B}$ (see Definition \ref{def:pwc_microstructures}), then, in the second step, using a density argument based on the approximation result Lemma \ref{lem:Lambda} to conclude in the general case of a jump set defining a Caccioppoli partition of $Q$, while in the last step follows from the density result of Caccioppoli partitions obtained in \cite[Lemma 5.3]{CriGra}.
As it usually happens, given a general $u\in L^2(\o; BV(Q;\{0,1\}))$, we cannot choose the sequence of piecewise constant functions that approximate it both in configuration and in energy. Instead, we need to construct it based on the function $u$ itself to dictate such piecewise approximation. This requires to have at our disposal a countable family $\mathcal{C}$ of sets of finite perimeter in $Q$ that are dense both in $L^1$ and also in energy.

The goal of Lemma \ref{lem:Lambda} is to construct such a countable family $\mathcal{C}$.
Note that we have to pay extra care in the construction of the recovery sequence, since the zeros $a, b:Q\to\R^M$ might be discontinuous on $\cup_{i=1}^k \partial E_i$.
Indeed, given $A\subset Q$ of finite perimeter, in order to approximate $A$ both in configuration and in energy we also need the approximating sequence $\{A_n\}_n\subset\mathcal{C}$ to be such that
\[
\lim_{n\to\infty}\hno\left( \left[\partial* A \triangle \partial* A_n\right] \cap
    \left(\,\bigcup_{i=1}^k \partial E_i\,\right) \right) =0.
\]
The countable family we choose to obtain the approximation is the following.

\begin{definition}
Let $\mathcal{C}$ be the family of open sets $G\subset Q$ such that there exist $\nu_1,\dots,\nu_m\in\S^{N-1}\cap \Q^M$, and $q_1,\dots,q_m\in\Q^M$, for which
\[
\partial G \subset \bigcup_{i=1}^m (q_i + \nu_i^\perp) \, \cup \, \bigcup_{i=1}^k \partial E_i,
\]
for some $m\in\N$.
\end{definition}

In the rest of the paper it will be convenient to adopt the following abuse of notation.
Given a set $A\subset Q$ with finite perimeter in $Q$, we will write $\widetilde{G}^1(A)$ in place of $\widetilde{G}^1(\widetilde{v_A})$, where $\widetilde{v}:Q\to\R^M$ is defined as
\[
\widetilde{v}_A(y)\coloneqq \chi_A(y)a(y) + [1-\chi_A(y)]b(y).
\]
Note that, using the fact that $A$ has finite perimeter, we have that $\widetilde{v_A}\in\widetilde{\mathcal{R}}$ (see Definition \ref{def:R}).

We are now in position to prove the first technical result.

\begin{lemma}\label{lem:Lambda}
For every $\lambda>0$ and every set of finite perimeter $A\subset Q$, there exists $E\in\mathcal{C}$ such that
\[
|A\triangle E| + \left| \widetilde{G}^1(A) - \widetilde{G}^1(E)  \right| < \lambda.
\]
\end{lemma}

\begin{proof}
\textbf{Step 1.} Without loss of generality, we can assume that
\begin{equation}\label{eq:zero_on_boundary}
\hno\left( \partial*A \cap \partial Q \right)
= \hno\left(\, \bigcup_{i=1}^k \partial E_i \cap \partial Q \,\right)= 0.
\end{equation}
Indeed, we can find $v\in\R^M$ such that the above condition is satisfied by the sets $A+v$ and $E_i +v$ in place of $A$ and $E_i$ respectively.
We then consider the sets $E_i+v$, $A+v$, and use the energy $\widetilde{G}_1^v$ defined as
\begin{equation}\label{eq:AF_close}
\widetilde{G}_1^v(B) \coloneqq \int_{\partial* B} \mathrm{d}_{\widetilde{W}}(y, \widetilde{v}^+(y), \widetilde{v}^-(y))\, \dhno(y),
\end{equation}
where $\widetilde{W}(y)\coloneqq W(y-v)$.\\

\textbf{Step 2.} Fix $\widetilde{\lambda}>0$, that will be chosen later.
Let $F\subset Q$ be the set of finite perimeter given by Theorem \ref{thm:strong_BV} relative to $A$ and $\widetilde{\lambda}$.
In particular,
\[
\| \ca_A - \ca_F \|_{L^1(Q)} < \widetilde{\lambda},
\]
and, using (iii), (v), and (vii) of Theorem \ref{thm:strong_BV}, we obtain
\begin{align*}
| \widetilde{G}_1(A) - \widetilde{G}_1(F) | &\leq
    \int_{Q\cap \partial*A\setminus \partial F} \dw(y,v_A^+(y), v_A^-(y))\, \dhno(y) \\
&\hspace{2cm}+\int_{Q\cap \partial F \setminus \partial*A} \dw(y,v_F^+(y), v_F^-(y))\, \dhno(y) \\
&\leq C_1\left[\, \hno\left( Q\cap \partial*A\setminus \partial F \right)
    + \hno\left( Q\cap \partial F \setminus \partial*A \right) \,\right] \\
&\leq C_1\left[\, |D\ca_A|( D\setminus C )
    + \hno\left( Q\cap \partial F \setminus \partial*A \right) \,\right] \\
&\leq 2 C_1\widetilde{\lambda},
\end{align*}
for some constant $C_1>0$ depending only on the wells $a$ and $b$.\\

\textbf{Step 3.} We now approximate the set $F$ with a set $G\in\mathcal{C}$.
Note that if the wells $a$ and $b$ were continuous, then the proof would be easier, since every piecewise-$C^1$ set in $Q$ can be approximated in the Hausdorff metric with a polyhedral set and, every polyhedral set in $Q$ can be approximated by a set in $\mathcal{C}$.
Due to the fact that the boundary of the approximated has to coincide as much as possible with the boundary of $F$ on $\cup_{i=1}^k \partial E_i$, the construction requires a more delicate argument.

First, we isolate the singularities of $\cup_{i=1}^k \partial E_i$ as follows.
It is possible to find $S_1,\dots,S_m\in\mathcal{C}$ with $\partial F$ orthogonal to $\partial S_i$ for each $i=1,\dots,m$, such that
\[
\bigcup_{i\neq j=1}^k \left(\, \partial E_i \cap \partial E_j  \,\right) \subset \bigcup_{i=1}^m S_i,
\]
and
\begin{equation}\label{eq:small_singular}
\sum_{i=1}^m \left[\, \hno( \partial S_i ) + \hno( \partial F \cap S_i )  \,\right] < \widetilde{\lambda}.
\end{equation}
Now we isolate the part of $\partial F$ on $\cup_{i=1}^k \partial E_i$. Set
\[
K \coloneqq \partial F \,\cap\, \left(\, \bigcup_{i=1}^k \partial E_i \,\right) \,\setminus\,
    \bigcup_{i=1}^m S_i.
\]
Recalling that $S_i\in\mathcal{C}$ for each $i=1,\dots,m$, and thus that each $S_i$ is open, we get that $K$ is compact.
By the outer regularity of the $\hno$ measure on $\cup_{i=1}^m \partial E_i$, it is then possible to find $R_1,\dots,R_n\in\mathcal{C}$ with $\partial F$ orthogonal to $\partial R_i$ for each $i=1,\dots,n$ such that
\[
K \subset \bigcup_{i=1}^n R_i \cap \bigcup_{i=1}^k \partial E_i,
\]
\begin{equation}\label{eq:small_all_Ri}
\sum_{i=1}^n \hno( \partial R_i ) < (1+\widetilde{\lambda}) \, 
    \hno \left(\, \bigcup_{i=1}^k \partial E_i \,\right),
\end{equation}
and
\begin{equation}\label{eq:small_R_i}
\sum_{i=1}^n \hno( \partial F \cap R_i )  < \widetilde{\lambda},
\end{equation}
and
\[
\hno\left(\, \bigcup_{i=1}^n R_i \cap \bigcup_{i=1}^k \partial E_i \, \setminus\, K \,\right) < \widetilde{\lambda}.
\]
Since $K$ is compact, there exist smooth open set $B\subset \cup_{i=1}^n R_i $ with $K\subset B$.
Without loss of generality, we can assume that $R_i\cap K\neq\emptyset$ for all $i\in\{1,\dots,n\}$.
Up to rearranging the order of the sets, we can assume that
\[
R_i \cap \widetilde{\partial} B \neq \emptyset
\]
if and only if $i\in\{1,\dots, n_1\}$, for some $n_1\leq n$, where $\widetilde{\partial} B$ is the relative boundary of $B$ in $\cup_{i=1}^k \partial E_i$.
In particular, this means that
\begin{equation}\label{eq;intersection_Ri}
\widetilde{\partial} K \subset \bigcup_{i=1}^{n_1} R_i.
\end{equation}
Thanks to \eqref{eq:small_all_Ri}, we have that
\begin{equation}\label{eq:small_selected_Ri}
\sum_{i=1}^{n_1}\hno( \partial R_i ) < C\widetilde{\lambda},
\end{equation}
for some $C>0$. 

To conclude, Let $T_1,\dots,T_p$ be the connected components of
\[
\partial F\setminus \left[\, \bigcup_{i=1}^{n_1}R_i \cup \bigcup_{i=1}^m S_i   \,\right].
\]
For each $i=1,\dots,p$, it is possible to find $\mu_i>0$ such that the sets $T_i + B(0,\mu_i)$ are pairwise disjoint. Consider the sets $F_1,\dots,F_p$, defined as
\[
F_i\coloneqq F\cap \left[\, T_i + B(0,\mu_i) \,\right].
\]
Find $E_i,\dots,E_p\in\mathcal{C}$ such that $E_i\subset T_i + B(0,\mu_i)$, and define
\[
E \coloneqq \left(\, F\setminus \bigcup_{i=1}^p \left[\, T_i + B(0,\mu_i) \,\right]  \,\right)
    \,\cup\, \bigcup_{i=1}^p E_i \,\cup\, \bigcup_{i=1}^m S_i \,\cup\, \bigcup_{i=1}^{n_1} R_i.
\]
Up to choosing $\mu_1,\dots,\mu_p>0$ and $\widetilde{\lambda}>0$ sufficiently small, thanks to \eqref{eq:small_singular}, \eqref{eq:small_R_i}, and \eqref{eq:small_selected_Ri}, we get
\[
|F\triangle E| + \left| \widetilde{G}^1(F) - \widetilde{G}^1(E)  \right| < C\widetilde{\lambda},
\]
for some $C>0$. Thus, by using step 1 and by selecting $\widetilde{\lambda}>$ sufficiently small, we conclude.
\end{proof}

We now introduce the class of piecewice constant functions for which we will construct the recovery sequence directly.

\begin{definition}\label{def:pwc_microstructures}
Denote by $\mathcal{B}$ the space of functions $u\in L^2(\o; L^1(Q;\R^M))$ such that
\[
u(x,y) = \sum_{i=1}^m u_i(y) \ca_{\o_i}(x),
\]
where, for each $i=1,\dots,m$, $\partial \o_i\cap \o$ is a polyhedral set, and
\[
u_i(y) \coloneqq \ca_{A(x)}a(y) + \left[\, 1 - \ca_{A(x)}(y)\,\right]b(y)
\]
for a set $A(x)\subset Q$ with finite perimeter in $Q$.
\end{definition}

We recall the main result that ensures the existence of a recovery sequence for a microstructure in the cube $Q$. The construction is based on \cite[Theorem 1.9]{CriGra}, with the additional complication of having to approximate the possible discontinuity of the wells on $\partial E_i\cap \partial E_j$ for $i\neq j$.

\begin{proposition}\label{prop:limsupQ}
Let $v\in L^1(Q;\R^M)$ be such that $v(y)\in\{a(y), b(y) \}$ for a.e. $y\in Q$, and such that $A\coloneqq \{u=a\}$ has finite perimeter in $Q$.
Then there exists a sequence of $Q$-periodic functions $\{v_n\}_n\subset W^{1,2}(Q;\R^M)$ with $v_n\to u$ in $L^1(Q;\R^M)$ such that
\[
\lim_{n\to\infty}\widetilde{G}_n(v_n) = \widetilde{G}^1(v),
\]
and $\sup_n \|v_n\|_{L^\infty}<\infty$.
\end{proposition}

\begin{remark}\label{rem:recovery_sequence}
Note that we can apply the mentioned result thanks to Theorem \ref{thm:boundEL}.
Moreover, the proof presented in \cite[Theorem 1.9]{CriGra} has to be adapted in order to take care of the fact that here we consider the perimeter in $Q$ seen as the periodic flat torus: this can be done as in step 1 of the proof of Lemma \ref{lem:Lambda}.\
Finally, considering zeros $a, b$ with possible discontinuities is not an issue, since it is simply possible to consider the functions $a_z,\dots,a_k$ and $b_1,\dots,b_k$ as separate zeros (multiple wells).
\end{remark}

We are now ready to prove that main result of this section.

\begin{proof}[Proof of Proposition \ref{prop:limsup}]
Without loss of generality, we can assume that $G^1(u)<\infty$,
otherwise there is nothing to prove.\\

\textbf{Step 1.} First assume $u\in\mathcal{B}$. Write it as
\[
u(x,y) = \sum_{i=1}^m u_i(y) \ca_{\o_i}(x),
\]
where, for each $i=1,\dots,m$, the set $\partial \o_i\cap \o$ is polyhedral, and
\[
u_i(y) \coloneqq \ca_{A_i}a(y) + \left[\, 1 - \ca_{A_i}(y)\,\right]b(y)
\]
for a set $A_i\in\mathcal{C}$.
For each $n\in\N$, consider a grid $\{Q_j^n\}_{j=1}^{k_n}$ of disjoint cubes of the form
\[
Q^n_j = \delta_n [0,1)^n + z_j,
\]
for some $z_j\in \delta_n\Z^N$, such that $Q^n_j\cap\o\neq\emptyset$.
For each $n\in\N$ and $i=1,\dots,m$, let
\[
\mathcal{I}^i_n \coloneqq \left\{\, j\in\{1,\dots,k_n\} \,:\, \overline{Q^n_j}\cap \partial \o_r = \emptyset, \text{ for all } r\neq i  \,\right\},
\]
and set
\[
\mathcal{I}_n \coloneqq \bigcup_{i=1}^m\mathcal{I}^i_n.
%\quad\quad\quad\mathcal{J}_n \coloneqq \{1,\dots, k_n \}\setminus \mathcal{I}_n.
\]
Since $\partial\o$ is regular, we have that
\begin{equation}\label{eq:asymptotics}
\lim_{n\to\infty} (\#\mathcal{I}_n)(\delta_n^N) = |\o|.
\end{equation}
Moreover, define
\[
\widetilde{\o}^n_i \coloneqq \bigcup_{j\in \mathcal{I}^i_n} Q^n_j,
\]
for each $i\in\{1,\dots, k_n \}$.
Let
\[
\widetilde{S}_n \coloneqq S_n\setminus\{ x\in\o \,:\, \mathrm{dist}(x,\partial S_n)\geq 3\e_n \}
\]
and, for each $n\in\N$, let $\varphi_n:\o\to[0,1]$ be such that
\[
\varphi_n\equiv1 \text{ on } \bigcup_{i=1}^m \widetilde{\o}_i, \quad\quad\quad\quad
\varphi^n_i\equiv0 \text{ on } \widetilde{S}_n,
\]
and satisfying
\begin{equation}\label{eq:estimate_gradient_cutoff}
|\nabla \varphi_n| \leq \frac{C}{\e_n},
\end{equation}
for some $C>0$. We remark that it is possible to construct such a family of cut-off functions satisfying this last estimate because
\[
\mathrm{dist}\left(\widetilde{\o}^n_i, \widetilde{\o}^n_j \right) \geq 2\delta_n
\]
whenever $i\neq j$.
For each $i\in\{1,\dots, m\}$, let $\{v^i_n\}_n\subset W^{1,2}(Q;\R^M)$ be the recovery sequence for the microstructure $u_i$ provided by Proposition \ref{prop:limsupQ}.
Note that each $v^i_n$ is $Q$-periodic, and that
\begin{equation}\label{eq:recovery_bounded}
\sup_n \|u^i_n\|_{L^\infty}<\infty.
\end{equation}
Define, for $n\in\N$, the function $u_n\in W^{1,2}(\o;\R^M)$ as
\[
u_n(x) \coloneqq \sum_{i=1}^m \left[\, \varphi_n(x)\ca_{\o_i}(x)v^i_n \left(\frac{x}{\delta_n}\right) + (1-\varphi_n(x))a\left(\frac{x}{\delta_n}\right) \,\right].
\]
Then it is easy to see that $u_n\to u$ in $L^1(\o)\times L^1(Q)$.
We now prove the convergence in energy.
Set
\[
S_n\coloneqq \o\setminus \bigcup_{i=1}^{k_n} \widetilde{\o}^n_i,
\]
and observe that
\begin{equation}\label{eq:Stilde_n}
|\widetilde{S}_n|\leq C\e_n,
\end{equation}
and that, thanks to \eqref{eq:estimate_gradient_cutoff}
\begin{equation}\label{eq:estimate_gradient_vn}
|\nabla u_n|\leq \frac{C}{\e_n}\quad\text{ in } S_n.
\end{equation}
Using \eqref{eq:recovery_bounded} and \eqref{eq:estimate_gradient_vn}, as $n\to\infty$ we get
\begin{align}\label{eq:energy_Sn}
&\int_{S_n} \left[\, \frac{\delta_n}{\e_n}
		 W\left( \frac{x}{\delta_n}, u_n(x) \right) 
		\,+ \e_n\delta_n|\nabla u_n(x)|^2 \,\right] \,dx \nonumber \\
&= \int_{\widetilde{S}_n} \left[\, \frac{\delta_n}{\e_n}
		W\left( \frac{x}{\delta_n}, u_n(x) \right)  + \e_n\delta_n|\nabla u_n(x)|^2\,\right] \,dx
+\int_{S_n\setminus \widetilde{S}_n} \e_n\delta_n\frac{|\nabla a|^2}{\delta_n^2} \,dx \nonumber \\
&\to0,
\end{align}
where, in the last step, the convergence of the first integral follows from \eqref{eq:Stilde_n}, while for the last integral from the fact that $a$ is Lipschitz, and $\frac{\e_n}{\delta_n}\to0$.
Thus, by using a change of variable and the fact that $\{v^i_n\}_n$, from \eqref{eq:energy_Sn} and \eqref{eq:asymptotics}, we get
\[
\lim_{n\to\infty} G_n(u_n) = G^1(u).
\]
This proves that $\{u_n\}_n$ is a recovery sequence.\\

\textbf{Step 2.} We now consider $u\in L^1(\o;BV(Q;\{a,b\}))$. We will construct a recovery sequence by using a diagonal argument. To be precise, fixed $\lambda>0$ we will construct a function $v^\lambda\in \mathcal{B}$ with
\[
\| v^\lambda - u \|_{L^1\times L^1} \leq C \lambda,\quad\quad\quad
\left|\, G^1(v^\lambda) - G^1(u)  \,\right| \leq C \lambda,
\]
for some constant $C>0$ independent of $\lambda$.
Thanks to step 1, we can find a sequence $\{v^\lambda_n\}_n\subset W^{1,2}(\o;\R^M)$ with $v^\lambda_n\to v^\lambda$ as $n\to\infty$ such that
\[
\lim_{n\to\infty}G_n(v^\lambda_n) = G^1(v^\lambda).
\]
The conclusion will follow by using the estimates above together with the arbitrariness of $\lambda$ and a diagonal argument.\\

We are now left with constructing the function $v^\lambda\in \mathcal{B}$.
Since $G^1(u)<\infty$, it is possible to find  $\mu>0$ such that
\begin{equation}\label{eq:small_sets_integral}
\int_E \widetilde{G}^1(u(x,\cdot))\, dx \leq \lambda,
\end{equation}
whenever $E\subset \o$ is a measurable set with $|E|\leq \mu$.
Without loss of generality, we can assume $\mu < \lambda$.
Let $\mathcal{C}=\{ F_i \}_{i\in\N}$ be the countable family given by Lemma \ref{lem:Lambda}.
The idea is to set the function $v^n\in\mathcal{B}$ as
\[
v^n(x,y) \coloneqq \sum_{i=1}^{n_0+1} v_{F_i}(y) \ca_{\o_i}(x),
\]
for some $n_0\in\N$, where
\[
v_{F_i}(y) \coloneqq \ca_{F_i}a(y) + \left[\, 1 - \ca_{F_i}(y)\,\right]b(y).
\]
The sets $\o_1,\dots\o_{i_n}$ will be defined in several steps.
For $i\in\N$, let
\begin{equation}\label{eq:wide_oi}
\widetilde{\o}_i \coloneqq \left\{\, x\in \o \,:\, |A(x)\triangle F_i|
    + |\widetilde{G}^1(A(x)) - \widetilde{G}^1(F_i)| \leq \frac{\lambda}{2^i} \,\right\}\, \setminus\, \bigcup_{j=1}^{i-1} \widetilde{\o}_j,
\end{equation}
where $A(x)\coloneqq \{ y\in Q \,:\, u(x,y)=a(y) \}$, and we set $\widetilde{\o}_{-1}\coloneqq \emptyset$.
Note that the sets $\widetilde{\o}_i$ are measurable, pairwise disjoint, and, thanks to Lemma \ref{lem:Lambda}, we also get that
\[
\o = \bigcup_{i=1}^\infty \widetilde{\o}_i.
\]
Let $n_0\in\N$ be such that
\[
\left|\, \o \, \setminus \, \bigcup_{i=1}^{n_0} \widetilde{\o}_i \, \right| \leq \frac{\mu}{2}.
\]
Let
\begin{equation}\label{eq:M}
M\coloneqq 1\, \lor\, \max_{i=1,\dots,n_0+1} \widetilde{G}^1(F_i).
\end{equation}
We claim that it is possible to construct a partition $\o_1,\dots,\o_{n_0+1}$ of polyhedral sets with
\begin{equation}\label{eq:approx_measure_oi}
|\widetilde{\o}_i \triangle \o_i| \leq \frac{\mu}{M}
\end{equation}
for all $i=1,\dots,n_0$, and such that
\begin{equation}\label{eq:approx_measure_on0}
|\o_{n_0+1}|\leq \frac{\mu}{M}.
\end{equation}
Indeed, since the sets $\widetilde{\o}_i$ are measurable, by the inner and the outer regularity of the Lebesgue measure, for each $i=1,\dots,n_0$, there exist a compact set $K_i $ and an open set $A_i$ with $K_i\subset \widetilde{\o}_i \subset A_i$ and
\[
|A_i\setminus K_i| \leq \frac{\mu}{2 M n_0}.
\]
By using smooth approximation of the characteristic function of $K_i$, we can find a polyhedral set $\o_i\subset A_i$ satisfying \eqref{eq:approx_measure_oi}. Starting from constructing $\o_1$ and for each $i=2,\dots,n_0$, substracting the union of the previous polyhedral sets from $\o_i$, we can assume that they are pairwise disjoint.
Finally, we define
\[
\o_{n_0+1}\coloneqq \o\setminus \bigcup_{i=1}^{n_0}\o_i.
\]
This partition of $\o$ satisfies all of the required properties.

We now show that the desired estimates hold. First of all, we note that
\begin{align*}
\| v^\lambda - u \|^2_{L^1\times L^1} &\leq \sum_{i=1}^{n_0}
    \int_{\o_i} \left[\,\int_{A(x)\triangle F_i} a(y)\, dy  
            + \int_{\left(Q\setminus A(x)\right)\triangle \left(Q\setminus F_i\right)} b(y)\, dy \, \right] \\
&\hspace{2cm} + \int_{\o_{n_0+1}} \| v^\lambda(x,\cdot) - u(x,\cdot) \|\, dx \\
&\leq C \sum_{i=1}^{n_0}\int_{\o_i} |A(x)\triangle F_i|\, dx + C |\o_{n_0+1}| \\
&\leq C \sum_{i=1}^{n_0} \left[\, \int_{\widetilde{\o}_i} |A(x)\triangle F_i|\, dx
    + C |\o_i\triangle \widetilde{\o}_i| \,\right] + C |\o_{n_0+1}| \\
&\leq C\lambda + C \mu\leq C\lambda
\end{align*}
where the previous to last inequality follows from \eqref{eq:wide_oi} and \eqref{eq:approx_measure_on0}, $C>0$ is a constant depending only the wells $a, b$ and on $|\Omega|$, and we recall that in the last inequality we used the fact that we are assuming, without loss of generality, that $\mu\leq\lambda$.
In a similar way, we have that
\begin{align*}
| G^1(v^\lambda) - G^1(u) | &\leq
    \sum_{i=1}^{n_0} \int_{\o_i} | \widetilde{G}^1(v^\lambda(x,\cdot))
        - \widetilde{G}^1(u(x,\cdot)) | \, dx \\
&\hspace{1cm}+ \int_{\o_{n_0+1}} | \widetilde{G}^1(v^\lambda(x,\cdot))
        - \widetilde{G}^1(u(x,\cdot)) | \, dx \\
&\leq \sum_{i=1}^{n_0} \int_{\widetilde{\o}_i} | \widetilde{G}^1(v^\lambda(x,\cdot))
        - \widetilde{G}^1(u(x,\cdot)) | \, dx \\
&\hspace{1cm}+ \int_{\o_i\triangle \widetilde{\o}_i} | \widetilde{G}^1(v^\lambda(x,\cdot))
        - \widetilde{G}^1(u(x,\cdot)) | \, dx \\
&\hspace{2cm}+ \int_{\o_{n_0+1}} | \widetilde{G}^1(v^\lambda(x,\cdot))
        - \widetilde{G}^1(u(x,\cdot)) | \, dx \\
&\leq C\lambda,
\end{align*}
where the last inequality follows from \eqref{eq:M}, \eqref{eq:approx_measure_oi}, and \eqref{eq:approx_measure_on0}, together with the fact that $\mu\leq \lambda$.
This provides the required function $v^\lambda$.\\

\textbf{Step 3.} We conclude for the case of a general case of $u\in \mathcal{R}$ by using step 2 together with \cite[Lemma 5.3]{CriGra}. This concludes the proof of the proposition.
\end{proof}

%%%%%%%%%%%%%%%%%%%%%%%%%%%%%%%%%%%%%%%%%%%%%%%%
%%%%%%%%%%%%%%%%%%%%%%%%%%%%%%%%%%%%%%%%%%%%%%%%
%%%%%%%%%%%%%%%%%%%%%%%%%%%%%%%%%%%%%%%%%%%%%%%%
%%%%%%%%%%%%%%%%%%%%%%%%%%%%%%%%%%%%%%%%%%%%%%%%

\section*{Acknowledgements}\
The research of Irene Fonseca and Likhit Ganedi was partially supported under grants NSF-DMS1411646 and NSF-DMS1906238.

\bibliographystyle{siam}
\bibliography{Bibliography}

\end{document}